  \documentclass[final,3p,times]{elsarticle}

\usepackage{amsmath, amssymb}
\usepackage{amsthm}
\usepackage{lmodern} 
\usepackage{algorithm}
\usepackage{algpseudocode}
\usepackage{epstopdf}
\usepackage{caption} 

\DeclareFontFamily{U}{futm}{}
\DeclareFontShape{U}{futm}{m}{n}{<-> fourier-bb}{}
\DeclareMathAlphabet{\fourierbb}{U}{futm}{m}{n}

\usepackage{booktabs}
\usepackage{tikz-cd}

\newcommand{\CX}{\mathcal{C}}
\newcommand{\DX}{\mathfrak{D}}
\newcommand{\abs}[1]{\left| #1 \right|}
\newcommand{\set}[1]{\left\{#1 \right\}}
\newcommand{\cbrc}[1]{\left( #1 \right)}
\newcommand{\brc}[1]{\left[ #1 \right]}
\newcommand{\norm}[1]{\left\lVert #1 \right\rVert}
\newcommand{\G}{\mathcal{G}}
\newcommand{\GA}{\Gamma}

\newcommand{\inc}[1]{ {#1}^{\mathrm{(inc)}}}
\newcommand{\bmsf}[1]{\boldsymbol{\mathsf{#1}}}
\newcommand{\rs}{\boldsymbol{r}}
\newcommand{\vare}{\varepsilon}
\newcommand{\nn}{\hat{\boldsymbol{n}}}
\newcommand{\rt}{\boldsymbol{r}'}
 
\newcommand{\Curl}{\nabla\times}
\newcommand{\Div}{\nabla\cdot}

\newcommand{\ed}{\mathbf{d}}
\newcommand{\Jv}{\mathbf{J}}
\newcommand{\Jf}{\boldsymbol{J}}
\newcommand{\Wf}[1]{\mathcal{W}^{(#1)}}
\newcommand{\Wv}[1]{\mathbf{W}^{(#1)}}

\newcommand{\pf}{\Phi}
\newcommand{\pd}{\overline{\mathbf{\Phi}} }
\newcommand{\Av}{\mathbf{A}}

\newcommand{\Ac}{\tilde{\mathbf{A}}}

\newcommand{\Ev}{\mathbf{E}}
\newcommand{\Ef}{\boldsymbol{E}}

\newcommand{\Hv}{\mathbf{H}}
\newcommand{\Hf}{\boldsymbol{H}}

 \newcommand{\Dv}{\mathbf{D}}
 \newcommand{\Df}{\boldsymbol{D}}
 
 \newcommand{\Bv}{\mathbf{B}}
\newcommand{\Bf}{\boldsymbol{B}}

\newcommand{\sL}{\mathcal{S}}
\newcommand{\dL}{\mathcal{D}}
\newcommand{\RR}{\abs{\rs-\rt}}

\newcommand{\vv}{\boldsymbol{v}}
\newcommand{\Par}{\partial}

\newcommand{\hd}{\fourierbb{H}}
\newcommand{\ded}{\fourierbb{D}}

\journal{Journal of Computational Physics}

\begin{document}

\begin{frontmatter}



\title{A Hybrid DEC-SIE Framework for Potential-Based Electromagnetic Analysis of Heterogeneous Media}



\author[label1]{Amgad Abdrabou\corref{cor1}} 
            
\cortext[cor1]{Corresponding author.}
\ead{aabdrab@purdue.edu}

\author[label1]{Luis J. Gomez} 
 
\affiliation[label1]{organization={Elmore Family School of Electrical and Computer Engineering, Purdue University},
         city={West Lafayette},
       postcode={47906},
            state={IN},
            country={USA}}

\begin{abstract}
 
Analyzing electromagnetic fields in complex, multi-material environments presents substantial computational challenges. To address these, we propose a hybrid numerical method that couples discrete exterior calculus (DEC) with surface integral equations (SIE) in the potential-based formulation of Maxwell's equations. The technique employs the magnetic vector and electric scalar potentials ($\Av$--$\Phi$) under the Lorenz gauge, offering natural compatibility with multi-physics couplings and inherent immunity to low-frequency breakdown. To effectively handle both bounded and unbounded regions, we divide the computational domain: the inhomogeneous interior is discretized using DEC, a coordinate-free framework that preserves topological invariants and enables structure-preserving discretization on unstructured meshes, while the homogeneous exterior is treated using SIEs, which inherently satisfy the radiation condition and eliminate the need for artificial domain truncation. A key contribution of this work is a scalar-component reformulation of the SIEs, which reduces the number of surface integral operators from fourteen to two by expressing the problem in terms of the Cartesian components of the vector potential and their normal derivatives. In the interior DEC domain, each component of $\Av$ is represented accordingly as a discrete 0-form. This is not a departure from the DEC framework, but rather an adaptation that mirrors established scalar-field treatments within DEC, preserves the underlying geometric structure, and aligns naturally with the scalar-component SIE representation at the interface. The result is a unified formulation in which the potentials remain differential-form quantities in the algebraic sense, yet are discretized component-wise for improved compatibility, numerical conditioning, and computational efficiency. The proposed hybrid method thus offers a physically consistent, structure-preserving, and efficient framework for solving electromagnetic scattering and radiation problems in complex geometries and heterogeneous materials, while avoiding the complexity of conventional vector-potential SIE formulations.

\end{abstract}




\end{frontmatter}

\section{Introduction}

The analysis of electromagnetic phenomena in heterogeneous and unbounded domains poses significant computational challenges. Differential equation (DE) and integral equation (IE) methods each offer distinct advantages, yet also face notable limitations. DE methods are particularly effective for modeling heterogeneous materials, as they directly solve Maxwell's equations in spatial form. However, they require artificial domain truncation with prescribed boundary conditions, complicating applications to open-region problems. Approximate truncation techniques such as absorbing boundary conditions (ABCs)~\cite{engquist1977} and perfectly matched layers (PMLs)~\cite{berenger1994} are widely used for this purpose. A more rigorous alternative couples DE solvers with IE methods, which represent the exterior domain via equivalent sources, yielding hybrid finite element–boundary integral (FE–BI) formulations that exactly truncate the finite element domain~\cite{sheng1998}.

A further challenge in DE solvers is the low-frequency breakdown~\cite{chew2014,li2016finite,liu2018potential}. As frequency approaches zero, direct discretizations of Maxwell’s equations become ill-conditioned due to the imbalance between electric and magnetic energy scales. Remedies include specialized basis functions (e.g., loop–tree or loop–star decompositions)~\cite{andersen2006hierarchical,graglia1997higher} or additional equations enforcing charge–current continuity~\cite{lee2005non,manges1995generalized}. More fundamental solutions reformulate the governing equations themselves~\cite{teixeira1999lattice,bochev2006principles,stern2020geometric}, often by introducing the magnetic vector potential ($\mathbf{A}$) and electric scalar potential ($\Phi$) formulation that is inherently immune to low-frequency breakdown. 

In this work, we adopt such a potential-based formulation within the framework of discrete exterior calculus (DEC). DEC provides a coordinate-free, geometrically intuitive framework for expressing Maxwell's equations in terms of differential forms, or \emph{$k$-forms}, which naturally integrate over $k$-dimensional geometric domains. These forms act as linear functionals on $k$-chains, and can be viewed as abstractions of physical measurement probes, producing scalar values through integration over lines, surfaces, or volumes rather than point-wise sampling~\cite{deschamps1981,flanders1989,warnick1997,frankel2011}. Crucially, DEC preserves key topological invariants and physical laws at the discrete level~\cite{hirani2003discrete,desbrun2005discrete}, and in the potential-based formulation, it avoids low-frequency breakdown~\cite{liu2018potential,boyuan2023}. Furthermore, DEC aligns naturally with multi-physics coupling and Hamiltonian descriptions of electromagnetism (e.g., the Aharonov–Bohm effect~\cite{aharonov59}). However, prior potential-based DEC implementations have relied exclusively on approximate domain truncation~\cite{boyuan2023}, which limits their applicability to open-region problems.

Coupling DEC with surface integral equation (SIE) truncation presents a major barrier: potential-based surface formulations for penetrable media are often highly complicated. Standard SIE methods may involve up to 14 distinct integral operators with varying singularity types~\cite{shanker2019,shanker2022}, resulting in substantial numerical and implementation burdens. Furthermore, in the DEC domain, the magnetic vector potential $\Av$ will be treated as a 1-form (following the vector representation of $\Av$ in the SIE), with components naturally associated with mesh edges.  While this conventional approach is mathematically elegant, this representation leads to considerable practical complexity at the DEC–SIE interface: matching the edge-based 1-form quantities to the boundary’s surface integral operators requires computing and coupling the large set of operator terms. For example, the conventional vector potential integral equation takes the form
\begin{equation}\label{eq:11}
    \Av(\rs) = \inc{\Av}(\rs) + \mathcal{S}[\hat{\boldsymbol{n}}(\rs')\times \nabla' \times\Av(\rs')]+ \nabla \times\mathcal{S}[\hat{\boldsymbol{n}}(\rs')\times \Av(\rs')]-\nabla\mathcal{S}[ \hat{\boldsymbol{n}}(\rs')\cdot \Av(\rs')]-\mathcal{S}[\nabla'\cdot \Av(\rs')],\quad  \rs \in \Omega_1
\end{equation}
where $\sL[u]$ denotes the single-layer potential~\cite{colton_kress98}:
\[
\sL[u](\rs) = \int\limits_{\Gamma} \mathcal{G}(\mathbf{r},\mathbf{r}') u(\mathbf{r}')\, d\Gamma(\mathbf{r}'),
\]
with
\[
\G(\rs,\rt)  = \frac{\exp{(i k_0 \RR)}}{4\pi \RR}.
\]
Here, $\Gamma$ is the boundary between interior and exterior domains, and $\hat{\boldsymbol{n}}(\rs')$ is the outward unit normal. The operator structure in~\eqref{eq:11} involves tangential, normal, divergence, and curl components of $\mathbf{A}$. 

 \textbf{In this work, we restrict attention to nonmagnetic materials with spatially homogeneous permeability; we show that this restriction admits a simplified operator structure.}  In \textbf{particular}, we propose a simplified representation treating each Cartesian component of the magnetic vector potential $A_\tau$ (\(\tau \in \{x,y,z\}\)) and its normal derivative on the surface as
\begin{equation}\label{eq:13}
    A_\tau(\rs) =  \inc{A}_\tau(\rs) - \mathcal{S}\!\left( \frac{\partial A_\tau(\rs')}{\partial \hat{\boldsymbol{n}}(\rs')}\right) + \mathcal{D}[A_\tau(\rs')],\quad  \rs \in \Omega_1
\end{equation}
where the double-layer potential is
\[
\dL[u](\rs)=  \int\limits_{\Gamma} \frac{\partial \G(\rs,\rt)}{\partial \nn(\rt)}\, u(\rt)\, d\GA(\rt).
\]
This mirrors the scalar potential equation
\begin{equation}\label{eq:14}
\Phi(\rs) = \inc{\Phi}(\rs) - \sL\!\left(\frac{\partial \Phi(\rs')}{\partial\nn(\rs')} \right) + \dL\!\left(\Phi(\rs')\right), \quad  \rs \in \Omega_1
\end{equation}
and introduces no additional operators beyond the classical single- and double-layer potentials. This reduction significantly decreases computational cost and system size while maintaining accuracy, as demonstrated in Section~\ref{sec:results}.

For the interior inhomogeneous region, we employ DEC, which, compared to traditional FEM, naturally supports unstructured meshes, avoids higher-order continuity requirements, and yields sparse, well-conditioned matrices~\cite{boyuan2023,rabina2014numerical}. Unlike standard hybrid DE–IE methods that enforce boundary continuity via tangential field components (e.g., \(\mathbf{n} \times \mathbf{E}\))~\cite{chew1998}, we use the potential-based equations~\eqref{eq:13}–\eqref{eq:14}. To ensure compatibility, we adopt an unconventional representation in DEC: instead of modeling $\mathbf{A}$ on the boundary as a 1-form, we treat its Cartesian components $A_x$, $A_y$, $A_z$ as independent 0-forms. This choice departs from the standard DEC formulation but offers significant advantages when coupling to SIEs. By aligning the DEC representation with the scalar component structure of the SIE on the truncation boundary, the method avoids the complexity of computing and matching a large set of surface integral operators---reducing the number from fourteen to two---while preserving accuracy. The approach maintains compatibility with the DEC framework, since each scalar component remains discretized over the primal mesh, yet it sidesteps the costly dual-mesh and edge-form coupling otherwise required for a 1-form treatment. The resulting formulation is both simpler and better conditioned, making it more suitable for large-scale electromagnetic scattering problems. Similar strategies appear in DEC treatments of Poisson and scalar wave problems, where scalar fields still benefit from the structure-preserving properties of the method.

In summary, this work unifies three key ideas:
\begin{enumerate}
    \item \textbf{Potential-based formulation}: A Lorenz-gauge representation of Maxwell’s equations in terms of $\mathbf{A}$ and $\Phi$, enabling stable, gauge-consistent modeling across broad frequency and material regimes. 
    \item \textbf{Hybrid DEC–SIE formulation}: An accurate radiation treatment in unbounded domains via SIE, coupled with DEC discretization that preserves geometric and topological structure.
    \item \textbf{Simplification of the SIE and DEC equations}: In contrast to established methods that require as many as 14 integral operators to enforce the necessary coupling conditions~\cite{shanker2019,shanker2022}, our formulation reduces this to only 2 integral operators. This simplification is enabled by representing $\mathbf{A}$ with scalar components and discretizing each as a 0-form in DEC for streamlined interface compatibility. 
\end{enumerate}

The remainder of this paper is organized as follows. Section~\ref{sec:theory} develops the continuous theoretical framework, beginning with the $\Av$–$\Phi$ formulation of Maxwell's equations and the associated operators. Section~\ref{sec:dec_sie} integrates DEC and SIE into a unified hybrid system. Section~\ref{sec:results} presents numerical validations in accuracy, convergence, and robustness. Section~\ref{sec:conclusions} concludes and outlines future directions. Additional derivations appear in the appendices: \ref{app1} reviews Maxwell’s equations in exterior calculus form, \ref{app2} constructs the Galerkin Hodge star, and \ref{app3} details the discrete potential SIEs and their method-of-moments implementation.

\section{Theoretical framework}\label{sec:theory}
In this section, we begin with a brief overview of Maxwell’s equations in the $\Av$-$\Phi$ formulation.  We then introduce a component-wise expression for the vector potential, formulated under the Lorenz gauge. Building on this formulation, we introduce the corresponding surface integral equations for the homogeneous exterior domain. For the inhomogeneous interior, a continuous exterior calculus representation that casts Maxwell's equations in the language of differential forms is presented. In this paper, the magnetic permeability is that of free space, and a time-harmonic dependence $\exp{(-i\omega t)}$  is assumed and suppressed.
\begin{figure}[ht!]
    \centering
    \includegraphics[width=0.25\linewidth]{ 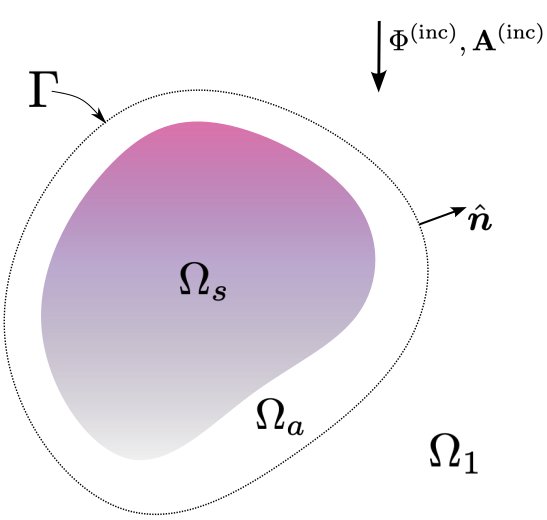}
    \caption{Geometry of a generic 3D object. The exterior domain $\Omega_1$ is free space, while the DEC domain $\Omega_2\equiv \Omega_s \cup \Omega_a$ is the union of two inhomogeneous domains $\Omega_s$ along with a buffer region $\Omega_a$. Due to the inclusion of the buffer region, the DEC domain boundary $\Gamma$ is completely in free space.}
    \label{Figure_1}
\end{figure}
We consider an inhomogeneous bounded region $ \Omega_s$ with spatially varying permitivity $\varepsilon(\rs) = \epsilon_0\varepsilon_s(\rs)$, where $\epsilon_0$ is the permittivity of free-space, and surrounded by free-space. Fields represented by a scalar potential $\Phi^{(inc)}$ and vector potential $\mathbf{A}^{(inc)}$ impinge on the object producing scalar potential $\Phi$ and vector potential $\mathbf{A}$ as shown in Fig.~\ref{Figure_1}. To find the potentials inside the inhomogeneous region, we define a computational domain $\Omega_2$ composed of the inhomogeneous bounded region $ \Omega_s$ and a homogeneous buffer $\Omega_a$ (i.e., $\Omega_2 \equiv \Omega_s \cup \Omega_a$) and the remaining homogeneous unbounded exterior region $\Omega_1$. By including the extra buffer region in the DEC domain $\Omega_2$, it is ensured that its boundary $\Gamma = \partial \Omega_2$ is embedded entirely in free space, which greatly simplifies the enforcement of boundary conditions.

\subsection{Maxwell's Equations in the $\Av$-$\Phi$ Formulation}
The electric and magnetic fields can be expressed in terms of the magnetic vector potential $\Av$ and electric scalar potential $\Phi$ as~\cite{harrington61}
\begin{equation}\label{eq:1}
    \Bv = \nabla \times \Av, \quad 
    \Ev = i\omega \Av - \nabla \Phi,
\end{equation}
where $\mathbf{E}$ and $\mathbf{B}$ denote electric field and magnetic flux density, respectively. The potentials then satisfy
\begin{align}
 \nabla \times \nabla \times \Av - \omega^2 \mu_0 \varepsilon_0 \varepsilon(\rs) \Av - i\omega \mu_0 \varepsilon_0 \varepsilon(\rs) \nabla \Phi &= \mu_0 \Jv, \label{eq:3}\\
 \nabla \cdot \left[\varepsilon(\rs) \nabla \Phi \right] - i \nabla \cdot \left[\omega \varepsilon(\rs) \Av \right] &= -\rho/\varepsilon_0, \label{eq:4}
\end{align}
where $\Jv$ and $\rho$ are the electric current and charge density, respectively. 
For notational and numerical convenience, we define the scaled magnetic vector potential $\Ac = \omega \Av$, with components $\Ac = \tilde{A}_x \hat{x} + \tilde{A}_y \hat{y} + \tilde{A}_z \hat{z}$. Using the continuity equation $\nabla \cdot \Jv = i\omega \rho$, with $k_0 = \omega\sqrt{\mu_0\varepsilon_0}$ and $\eta_0 = \sqrt{\mu_0/\varepsilon_0}$, we rewrite Eqs.~\eqref{eq:3} and~\eqref{eq:4} as:
\begin{align}
\nabla \times \nabla \times \Ac - k_0^2 \varepsilon(\rs) \Ac - i k_0^2 \varepsilon(\rs) \nabla \Phi &= k_0 \eta_0 \Jv, \label{eq:A}\\
\nabla \cdot \left[ \varepsilon(\rs) \nabla \Phi \right] - i \nabla \cdot \left[ \varepsilon(\rs) \Ac \right] &= \frac{i \eta_0}{k_0} \nabla \cdot \Jv.
\end{align}
Upon imposing the Lorenz gauge $\nabla \cdot \Ac = i k_0^2 \varepsilon(\rs) \Phi$, the system becomes:
\begin{equation}\label{eq:Contin}
   \begin{bmatrix}
        \mathcal{V} & 0 & 0 & -i\partial_x \varepsilon(\rs) \\
         0 & \mathcal{V} & 0 & -i\partial_y \varepsilon(\rs) \\
        0 & 0 & \mathcal{V} & -i\partial_z \varepsilon(\rs) \\
        -ik_0^2\partial_x \varepsilon(\rs) & -ik_0^2\partial_y \varepsilon(\rs) & -ik_0^2\partial_z \varepsilon(\rs) & \mathcal{L}
   \end{bmatrix}
   \begin{bmatrix}
       \tilde{A}_x \\ \tilde{A}_y \\ \tilde{A}_z \\ \Phi_s
   \end{bmatrix}
   = -\eta_0 k_0
   \begin{bmatrix}
       J_x \\ J_y \\ J_z \\ \frac{\nabla \cdot \Jv}{i}
   \end{bmatrix},
\end{equation}
where $\Phi_s := k_0^2 \Phi$, and the operators are
\begin{align}
\mathcal{V} \tilde{A}_\tau &:= \nabla^2 \tilde{A}_\tau + k_0^2 \varepsilon(\rs) \tilde{A}_\tau, \label{eq:OPV}\\
\mathcal{L} \Phi_s &:= \nabla \cdot \left[ \varepsilon(\rs) \nabla \Phi_s \right] + k_0^2 \varepsilon^2(\rs) \Phi_s. \label{eq:OPL}
\end{align}

Before introducing the DEC formulation of Eq.~\eqref{eq:Contin} and how the SIE is incorporated, we first need to cast Eq.~\eqref{eq:Contin} into the exterior calculus notation~\cite{hirani2003discrete}. This is achieved by first recognizing that differential 0-forms are equivalent to scalar fields. Hence, $\tilde{A}_x,\,\tilde{A}_y,\,\tilde{A}_z$, and $ {\Phi}_s$ can be viewed as 0-forms. Second, only two operators from exterior calculus are needed to provide the exterior calculus representation of many vector calculus operators, namely, the exterior derivative $\ed$ and the Hodge star $\star$ operators. The operator $\ed$ maps a $k$-form to a $(k+1)$-form, while the Hodge star enables us to establish the constitutive relations in electromagnetics, $\Dv = \vare_0 \vare(r) \Ev, \Bv = \mu_0\ \mu(r) \Hv$ as $\Df= \vare_0\star_\epsilon \Ef$ and $\Bf = \mu_0\star_\mu \Hf$, where $\Ef$  and $\Hf$ are 1-forms and $\Df$ and $\Bf$ are 2-forms~\cite{deschamps1981}. Last, for a vector field $\Av$ in $\mathbb{R}^3$ and a scalar field $\Phi$, the following transformations hold~\cite{mohammed2016}, which provide a way to bridge the 1-form $\boldsymbol{A}$ to vector field $\Av$ and the 0-form $\Phi$ with is also a scalar field as $\Curl \Av   \longrightarrow \star\ed\boldsymbol{A}$, $\Div \Av \longrightarrow \star\ed\star\boldsymbol{A}$, and $\nabla \Phi\longrightarrow \ed \Phi$. Using such a framework, we can now write the system~\eqref{eq:Contin} in the exterior calculus notation as 
\begin{equation}
    \label{eq:Contin_DEC}
   \renewcommand{\arraystretch}{1.2} 
    \setlength{\arraycolsep}{6pt} 
    \begin{bmatrix}
        \mathcal{V}_{\mathrm{EC}}              & 0  & 0 & -i\star_{\partial_x \vare(\rs)}\\
         0  & \mathcal{V}_{\mathrm{EC}}  &  0 & -i\star_{\partial_y \vare(\rs)} \\
         0  &  0  & \mathcal{V}_{\mathrm{EC}}  &  -i\star_{\partial_z \vare(\rs)}\\
        -ik_0^2\star_{\partial_x \vare(\rs)} &  -ik_0^2\star_{\partial_y \vare(\rs)}  & -ik_0^2\star_{\partial_z \vare(\rs)}& \mathcal{L}_{\mathrm{EC}} \\
    \end{bmatrix}  \begin{bmatrix}
        \tilde{A}_x\\
       \tilde{A}_y\\
      \tilde{A}_z\\
       {\Phi}_s
    \end{bmatrix} = -\eta_0 k_0 \star \begin{bmatrix}
        J_x\\
         J_y\\
        J_z\\
        \frac{\Div \Jv}{i}
    \end{bmatrix},
\end{equation}
with
\begin{equation}\label{eq:OPVd}
\mathcal{V}_{\mathrm{EC}}  :=  \ed \star \ed + k_0^2 \star_{\vare} ,
\end{equation}
 and
\begin{equation}\label{eq:OPLd}
\mathcal{L}_{\mathrm{EC}}    := \ed \star_{\vare} \ed + k_0^2 \star_{\vare^2}
\end{equation}
\ref{app1} offers a brief introduction to the exterior calculus notation and the definition of exterior derivative and Hodge star operators.
\subsection{Surface integral equations in the exterior region}\label{sec2_2}
The surface integral equation decouples the exterior domain $\Omega_1$ from the interior one. This is done by introducing equivalent sources. Furthermore, since the exterior domain is homogeneous, the fields can be separated into incident fields $\tilde{\Av}^{(\mathrm{inc})}$ and $\Phi^{(\mathrm{inc})}$ and scattered fields generated by the equivalent sources. The incident fields are either due to impressed sources $\Jv$ and $\rho$ as
\[
 \tilde{\Av}^{(\mathrm{inc})}(\rs) = k_0\eta_0 \mathcal{S}[\Jv](\rs), \quad 
 \Phi^{(\mathrm{inc})}(\rs) = \frac{1}{\varepsilon_0} \mathcal{S}[\rho](\rs).
\]
or for an incident plane wave~\cite{vico2016}:
\begin{align}\label{eq:incAPhi}
    \Phi^{(\mathrm{inc})}(\rs) &= -(\rs \cdot \Ev_p)\, e^{i k_0 \hat{\boldsymbol{u}} \cdot \rs}, \\
    \tilde{\Av}^{(\mathrm{inc})}(\rs) &= -k_0 (\rs \cdot \Ev_p)\, e^{i k_0 \hat{\boldsymbol{u}} \cdot \rs}\, \hat{\boldsymbol{u}}.
\end{align}
This representation of sources is compatible with field equations (i.e., $\Ev^{\mathrm{inc}} = -\nabla \Phi^{\mathrm{inc}} + i \tilde{\Av}^{(\mathrm{inc})}$). The jump conditions from~\cite{colton_kress98}, for each scalar component, are the following:
\begin{align}
       \frac{1}{2}\tilde{A}_x(\rs) - \widetilde{\mathcal{D}}[\tilde{A}_x](\rs) + \mathcal{S}[\frac{\partial } {\partial n} \tilde{A}_x](\rs)& =  \tilde{A}_x^{(\mathrm{inc})}(\rs),  \label{eq:SIE_Ax}\\
       \frac{1}{2}\tilde{A}_y(\rs) - \widetilde{\mathcal{D}}[\tilde{A}_y](\rs) + \mathcal{S}[\frac{\partial } {\partial n} \tilde{A}_y](\rs)& =  \tilde{A}_y^{(\mathrm{inc})}(\rs),  \label{eq:SIE_Ay}\\
       \frac{1}{2}\tilde{A}_z(\rs) - \widetilde{\mathcal{D}}[\tilde{A}_z](\rs) + \mathcal{S}[\frac{\partial } {\partial n} \tilde{A}_z](\rs)& =  \tilde{A}_z^{(\mathrm{inc})}(\rs),  \label{eq:SIE_Az}\\
       \frac{1}{2}\Phi(\rs) - \widetilde{\mathcal{D}}[\Phi](\rs) + \mathcal{S}[\frac{\partial } {\partial n} \Phi](\rs)& =  \Phi^{(\mathrm{inc})}(\rs),  \label{eq:SIE_Phi}
\end{align}

where $ \rs \in \Gamma$. Only the classical single-layer $\mathcal{S}$ and double-layer $\widetilde{\mathcal{D}}$ operators are required, greatly simplifying the surface formulation.  Furthermore, since the surface $\Gamma$ is entirely embedded in free-space, the boundary conditions for the electric scalar potential and each component of the vector potential are naturally continuous across $\Gamma$. 
\section{The hybrid DEC-SIE formulation}\label{sec:dec_sie}
 We now introduce the hybrid DEC-SIE formulation, combining exterior calculus for the inhomogeneous interior domain with a simplified SIE approach for the homogeneous exterior. We start with the DEC framework that is used to discretize the fields in the interior domain, then briefly present the discrete SIEs. Finally, the hybrid DEC-SIE approach is presented.  
 \subsection{The DEC framework}\label{subsec:dec}
 Unlike finite element methods, which require the explicit selection of basis and testing functions, DEC implicitly defines the approximation of electromagnetic quantities through the choice of a simplicial complex that discretizes the geometry. Like with finite element methods, for the implementation of DEC the interior domain $\Omega_2$ (see Fig.~\ref{Figure_1}) is first discretized using a tetrahedral mesh known as a simplicial complex, denoted by $\mathcal{K}$. The complex $\mathcal{K}$ contains a natural hierarchy of simplices: 3-simplices (tetrahedra), their 2-simplex faces (triangles), their 1-simplex edges, and their 0-simplex nodes. Each $k$-simplex $\sigma^k \in \mathcal{K}$, for $k = 0, \ldots, 3$, is defined by its $k+1$ vertices as $\sigma^k = [v_1, \ldots, v_{k+1}]$, where the subscripts denote node indices \cite{hirani2003discrete}. The 0-simplices are the vertices $\{v_1, \ldots, v_{N_0}\}$, where $N_0$ is the number of nodes. The 1-simplices are edges $\{e_1, \ldots, e_{N_1}\}$, with each edge $e_k \equiv \sigma^1_k = [v_{k_1}, v_{k_2}]$ oriented according to the ordering of its endpoints, and $N_1$ is the number of edges. The 2-simplices are triangles $\{t_1, \ldots, t_{N_2}\}$, with each triangle $t_k \equiv \sigma^2_k = [v_{k_1}, v_{k_2}, v_{k_3}]$ defined by three nodes and bounded by three 1-simplex edges. The 3-simplices are tetrahedra $\{\mathcal{T}_1, \ldots, \mathcal{T}_{N_3}\}$, where each $\mathcal{T}_k \equiv \sigma^3_k = [v_{k_1}, v_{k_2}, v_{k_3}, v_{k_4}]$ consists of four nodes and is bounded by four triangular faces (2-simplices). The total number of tetrahedra is denoted $N_3$. The orientation of the 3-simplices is assumed to be consistent throughout the mesh, a condition typically enforced by mesh generation tools. For lower-dimensional simplices ($k < 3$), we adopt a canonical orientation: each edge $[v_{k_1}, v_{k_2}]$ is oriented such that $k_1 < k_2$, and each triangle $[v_{k_1}, v_{k_2}, v_{k_3}]$ such that $k_1 < k_2 < k_3$. 

Since the simulation domain is now represented by a simplicial complex $\mathcal{K}$, the smooth differential forms must also be discretized. In the context of DEC, this discretization corresponds to evaluating the integral of each smooth differential $k$-form over the $k$-simplices of the mesh. 
The corresponding spaces of discrete differential forms are denoted by $\mathcal{C}^k(\mathcal{K})$ for $k$-forms defined on $\mathcal{K}$. Discretization of smooth differential forms is carried out via the de Rham map $\mathcal{R}$, which integrates the form over the appropriate $k$-simplices. For instance, the maps that discretize the primal 0-form $\pf$, 1-form $\Ef$, and 2-form $\Bf$ into their corresponding degree-of-freedom arrays are given by:
\begin{equation}\label{eq:deRhamMapsPrime}
\begin{aligned}
\mathcal{R}(\pf) = \overline{\boldsymbol{\Phi}} = \left[ \Phi(v_1), \ldots, \Phi(v_{N_0}) \right]^\mathsf{T} \in  \mathcal{C}^0(\mathcal{K}), \\
\mathcal{R}(\Ef)= \overline{\mathbf{e}} = \left[ \int_{e_1} \Ef, \ldots, \int_{e_{N_1}} \Ef \right]^\mathsf{T} \in \mathcal{C}^1(\mathcal{K}), \\
\mathcal{R}(\Bf)  = \overline{\mathbf{b}} = \left[ \int_{t_1} \Bf, \ldots, \int_{t_{N_2}} \Bf \right]^\mathsf{T}\in\mathcal{C}^2(\mathcal{K}),
\end{aligned}
\end{equation}
where the superscript $\mathsf{T}$ denotes the transpose. 

Next, we provide DEC representations of the exterior derivative $\ed$ and $\star$ operators. The discrete exterior derivative operator, denoted by $\ded_k$, maps discrete $k$-forms to $(k+1)$-forms on the complex $\mathcal{K}$. It encodes the combinatorial structure of the mesh and is derived directly from the generalized Stokes' theorem
\begin{equation}\label{eq:stokes}
  \int_{\partial \mathcal{M}} \boldsymbol{\alpha} = \int_{\mathcal{M}} \ed \boldsymbol{\alpha},
\end{equation}
where $\alpha$ is a $k$-form, $\mathcal{M}$ is a $(k+1)$-dimensional oriented manifold, and $\partial \mathcal{M}$ is its boundary. Thus, the degree $k$ of a differential $k$-form determines the dimensionality of the geometric object over which it is integrated. For electromagnetic fields, the manifold $\mathcal{M} \subseteq \mathbb{R}^3$, the maximum degree of forms is limited by $k \leq 3$. To illustrate, consider the smooth 0-form $\pf$ and its exterior derivative $\ed\pf$, which is a 1-form. When discretized, $\pf$ is evaluated at the nodes, while $\ed\pf$ is integrated along each oriented edge $e_k = [v_{k_1}, v_{k_2}]$. By Stokes’s theorem, we have:
\[
[\mathcal{R}(\ed\pf)]_k = \int_{e_k} \ed\pf = \pf(v_{k_2}) - \pf(v_{k_1}) = \sum_\ell [\ded_0]_{k\ell} \Phi_\ell,
\]
where $\Phi_\ell = \pf(v_\ell)$ and $k = 1,\ldots,N_1$. The matrix $\ded_0$ is an $N_1 \times N_0$ incidence matrix defined as
\[
[\ded_0]_{k\ell} =
\begin{cases}
+1, & \text{if node $\ell$ is the head of edge $k$}, \\
-1, & \text{if node $\ell$ is the tail of edge $k$}, \\
0,  & \text{otherwise}.
\end{cases}
\]
A similar construction applies to higher-degree forms, where  $\ded_1$ acting on discrete 1-forms is an $N_2 \times N_1$ matrix defined as
\[
[\ded_1]_{k\ell} =
\begin{cases}
+1, & \text{if edge $\ell$ is on the boundary of triangle $k$ with matching orientation}, \\
-1, & \text{if edge $\ell$ is on the boundary of triangle $k$ with opposite orientation}, \\
0,  & \text{otherwise}.
\end{cases}
\]
and $\ded_2$ acting on discrete 2-forms is an $N_3 \times N_2$ matrix defined as
\[
[\ded_2]_{k\ell} =
\begin{cases}
+1, & \text{if triangle $\ell$ is a face of tetrahedron $k$ with consistent orientation}, \\
-1, & \text{if triangle $\ell$ is a face of tetrahedron $k$ with opposite orientation}, \\
0,  & \text{otherwise}.
\end{cases}
\]

Certain operations---notably the Hodge star, which maps $k$-forms to $(n-k)$-forms, where $n$ is the dimension of $\mathcal{K}$---require an additional geometric structure: the dual complex $\ast\mathcal{K}$. For $n=3$, the dual complex $\ast\mathcal{K}$ assigns to each $k$-simplex in $\mathcal{K}$ a dual $(3-k)$-cell, constructed so that the primal and dual meshes are orthogonal. This construction allows integration over volumes, surfaces, or paths to be transferred between primal and dual complexes. Alternatively, Whitney-form-based Hodge stars define the mapping directly via a variational formulation without explicitly requiring a dual mesh \cite{Bossavit1999}. Barycentric-based Whitney approaches can be used to construct sparse, geometry-dependent Hodge matrices \cite{Kurz}. In this work, we adopt the Galerkin formulation for the Hodge star, which ensures consistency and flexibility without explicitly constructing $\ast\mathcal{K}$. The Galerkin Hodge star operator is denoted by $\hd_k$, and maps a discrete differential $k$-form to a $(3-k)$-form. For the constitutive relations $\Df = \vare_0\star_\vare \Ef$ and $\Hf ={\mu_0}^{-1}\star_{\mu^{-1}}\Bf$, we designate the material parameters dependence as
\[
\overline{\mathbf{d}} = \vare_0\hd_1(\vare) \overline{\mathbf{e}}, \quad \overline{\mathbf{h}} = \mu_0^{-1}\hd_2(\mu^{-1}) \overline{\mathbf{b}},
\]
where $\hd_k(\cdot)$ is an $N_k\times N_k$ matrix, $k = 0,\ldots,3$. The construction of the matrices $\hd_k$ involves the use of Whitney forms $\Wf{k}(\rs)$, or their proxy fields obtained via the sharp $\sharp$ operator that maps 1-forms to vector fields. The definitions of such forms and the process of constructing $\hd_k$ are explained in~\ref{app2}. For details on Whitney forms and their proxy fields, the reader may refer to~\cite{whiney1957,lohi2021}. The orientation of each primal $k$-simplex induces a compatible orientation on the corresponding dual $(3{-}k)$-cell~\cite{desbrun2005discrete}. As a result, the discrete exterior derivative on the dual mesh, denoted by $\ded_k^\ast$, is related to the transpose of the primal operator as $\ded_k^\ast = (-1)^{3-k} \ded_{3-k-1}^\mathsf{T}.$
These relations enable computations on the dual mesh without explicitly constructing it, and are needed to approximate the operators $\mathcal{V}_{\mathrm{EC}}$ and $\mathcal{L}_{\mathrm{EC}}$ defined by Eqs.~\eqref{eq:OPVd} and~\eqref{eq:OPLd}.

Although our formulation only involves 0-forms as unknowns, it is important to understand how discrete differential forms of higher degree arise. For example, the exterior derivative $\ed\Phi$ is a 1-form and is discretized on edges, while its Hodge dual $\star\Phi$ is a 3-form and is associated with volumes in the dual mesh (corresponding to nodes in the primal mesh). The following diagram shows how the discrete exterior derivative and Hodge star operators relate the spaces of discrete differential forms for the 3D case:
\begin{equation}\label{eq:diagram}
\begin{tikzcd}[column sep=large, row sep=large ]
 \CX^0(\mathcal{K}) \arrow{r}{\ded_0} \arrow[d, "\hd_0", shift left=1ex ] & \CX^1(\mathcal{K}) \arrow{r}{\ded_1} \arrow[d, "\hd_1", shift left=1ex ]  & \CX^2(\mathcal{K}) \arrow{r}{\ded_2} \arrow[d, "\hd_2", shift left=1ex ] & \CX^3(\mathcal{K})    \arrow[d, "\hd_3", shift left=1ex ]\\
 \DX^3(\mathbf{\ast} \mathcal{K})  \arrow[u, "\hd_0^{-1}", shift left=1ex ] &  \DX^2(\ast \mathcal{K}) \arrow[l, "-\ded_0^\mathsf{T}"{swap} ] \arrow[u, "\hd_1^{-1}", shift left=1ex ]  &  \DX^1(\ast \mathcal{K}) \arrow[l, "\ded_1^\mathsf{T}"{swap} ] \arrow[u, "\hd_2^{-1}", shift left=1ex ] &  \DX^0(\ast \mathcal{K})   \arrow[l, "-\ded_2^\mathsf{T}"{swap} ] \arrow[u, "\hd_3^{-1}", shift left=1ex ]
\end{tikzcd}
 \end{equation}
 where $\mathfrak{D}^k(\ast\mathcal{K})$ is the space of discrete differential $k$-forms defined on $\ast\mathcal{K}$. To this point, we have introduced all the necessary discrete exterior calculus representations of the exterior calculus operators $\ed$ and $\star$. 
\subsection{The discrete potential SIEs}\label{subsec:sie}
On the boundary of the $\mathcal{K}$ that approximates $\Gamma$, the system of integral equations~\eqref{eq:SIE_Ax}-\eqref{eq:SIE_Phi} is discretized using the method of moments (MoMs)~\cite{harrington68}. The boundary mesh consists of $N_2^\Par$ boundary triangles, $N_1^\Par$ boundary edges, and $N_0^\Par$ boundary nodes. Applying the MoMs yields the following system of equations:
\begin{align}
  \bmsf{D}\,{\mathbf{a}}_x^{(\Gamma)}+ \bmsf{S}\,\partial_{\nn}{\mathbf{a}}_x^{(\Gamma)} &=   \boldsymbol{f}_{\tilde{A_x}}, \label{eq:main_Ax}\\
   \bmsf{D}\,{\mathbf{a}}_y^{(\Gamma)}+ \bmsf{S}\,\partial_{\nn}{\mathbf{a}}_y^{(\Gamma)} &= \boldsymbol{f}_{\tilde{A_y}}, \label{eq:main_Ay}\\
   \bmsf{D}\,{\mathbf{a}}_z^{(\Gamma)}+ \bmsf{S}\,\partial_{\nn}{\mathbf{a}}_z^{(\Gamma)} &= \boldsymbol{f}_{\tilde{A_z}}, \label{eq:main_Az}\\
      \bmsf{D}\, {\mathbf{\Phi}}^{(\Gamma)}+ \bmsf{S}\,\partial_{\nn}{\mathbf{\Phi}}^{(\Gamma)} &= \boldsymbol{f}_{\Phi}, \label{eq:main_Phi}
\end{align}
where $\mathbf{a}_x^{(\Gamma)},\mathbf{a}_y^{(\Gamma)},\mathbf{a}_z^{(\Gamma)}$ and $\Phi^{(\Gamma)}$ are $N_0^\partial$-dimensional column vectors presenting the evaluations of the scalar fields $\tilde{A}_x(\rs_n), \tilde{A}_y(\rs_n), \tilde{A}_z(\rs_n)$, and $\Phi(\rs)$, respectively, with $\rs_n$, $n = 1,\ldots, N_0^\Par$ denoting position vectors of the boundary nodes. The respective normal derivatives of these fields at the boundary nodes are denoted by $\partial_{\nn}{\mathbf{a}}_x^{(\Gamma)}$, $\partial_{\nn}{\mathbf{a}}_x^{(\Gamma)}$, $\partial_{\nn}{\mathbf{a}}_x^{(\Gamma)}$, and $\partial_{\nn}{\mathbf{\Phi}}_s^{(\Gamma)}$. The $N_0^\Par\times N_0^\Par$ matrices $\bmsf{D}$ and $\bmsf{S}$ are the MoM approximations to the double and single-layer potentials, and their entries, along with the $N_0^\partial$-dimensional column vectors $\boldsymbol{f}_{\tilde{A}_{\tau}}, \tau\in\set{x,y,z}$ and $\boldsymbol{f}_{\Phi}$ are given in~\ref{app3}.
\subsection{The linear system}\label{subsec:hybrid}
 To obtain the discrete representation of the governing system~\eqref{eq:Contin_DEC} within the DEC formulation and incorporate the SIEs~\eqref{eq:main_Ax}-\eqref{eq:main_Phi} into the DEC framework, we proceed as follows. The continuous 0-forms $\tilde{A}_\tau$, $\tau \in \{x, y, z\}$, and $\Phi_s$ are approximated by discrete 0-forms $\overline{\mathbf{a}}_\tau$ and $\overline{\mathbf{\Phi}}_s$, respectively. These discrete variables represent the evaluations of the corresponding continuous 0-forms at the nodes of the primal mesh, as discussed in Section~\ref{subsec:dec}. With these approximations, the system of equations~\eqref{eq:Contin_DEC} can be fully expressed within the DEC framework. For simplicity, we assume that the excitation arises solely from the incident fields in $\Omega_1$ and set $\mathbf{J} = 0$. Accordingly, in the DEC domain, we obtain
\begin{equation}\label{eq:Atau}
\brc{ -\ded_0^{\mathsf{T}}\hd_1 \ded_0 + k_0^2\hd_0(\vare)  }\overline{\mathbf{a}}_\tau + \ded_2^{\mathrm{(c)}}\partial_{\nn}{\mathbf{a}}_\tau^{(\Gamma)}  - i \hd_0\cbrc{\Par_\tau\vare }\overline{\mathbf{\Phi}}_s = 0,\quad \forall\,\tau\in\set{x,y,z}
\end{equation}
and
\begin{equation}\label{eq:Phitau}
\brc{-\ded_0^{\mathsf{T}}\hd_1(\vare) \ded_0 + k_0^2\hd_0(\vare^2) } \overline{\mathbf{\Phi}}_s +\ded_2^{\mathrm{(c)}}\partial_{\nn}{\mathbf{\Phi}}_s^{(\Gamma)}- ik_0^2\sum\limits_{\tau} \hd_0\cbrc{\Par_\tau \vare}\overline{\mathbf{a}}_\tau = 0.
\end{equation}
Here, the terms $\ded_2^{\mathrm{(c)}}\partial_{\nn}{\mathbf{a}}_\tau^{(\Gamma)}$ and $\ded_2^{\mathrm{(c)}}\partial_{\nn}{\mathbf{\Phi}}_s^{(\Gamma)} $ arise due to the incomplete dual mesh at the boundary $\Gamma$, and are essential for coupling the DEC formulation with surface integral equations. The operator $\hd_0\cbrc{\Par_\tau \vare}$ is constructed using the Galerkin approach with Whitney 0-forms and reduces to surface integrals on the common faces of the tetrahedra. Assuming $\Gamma$ is approximated by a surface triangulation with $N_0^\partial$ boundary nodes number, the operator $\ded_2^{\mathrm{(c)}}$ complements the dual exterior derivative operator $-\ded_0^{\mathsf{T}}$. To define the complement operator $\ded_2^{\mathrm{(c)}}$, we first observe that the operator $-\ded_0^{\mathsf{T}}$ acts on dual 2-cochains~\footnote{The notions of cochain and discrete differential form are equivalent and used interchangeably.}. For example, in the term
\[
 -\ded_0^{\mathsf{T}}  \hd_1\ded_0\pd_s,
\]
the quantity $\hd_1\ded_0\pd_s$ is a dual 2-cochain. However, due to the incompleteness of the dual mesh near the boundary, the operator $-\ded_0^{\mathsf{T}}$ alone is insufficient; it must be complemented with boundary contributions~\cite{mohammed2016}.  To motivate the required complement, we invoke the identity $\star \ed \Phi = (\nabla\Phi \cdot \boldsymbol{n})\nu$, where $\nu$ is the surface 2-form on the boundary $\Gamma$~\cite[p.~506]{marsden2012}. This suggests that boundary flux terms supplement the incomplete dual structure. Accordingly, the operator $\ded_2^{\mathrm{(c)}}$ is an $N_0\times N_0^\partial$ matrix that augments the dual 3-cells by traversing along the primal boundary triangles in the orientation direction of the dual 3-cells, and is given by
\begin{equation}\label{eq:d2C}
 \ded_2^{\mathrm{(c)}} = \frac{1}{6} \, \abs{\ded_0^\mathsf{T}} \, \abs{\ded_1^\mathsf{T}} \, \mathbf{P}_2^\mathsf{T} \, \Pi_S \, \mathbb{Q}, 
\end{equation}
where $\abs{\ded_k^\mathsf{T}}$ denotes the matrix $\ded_k^\mathsf{T}$ with all entries replaced by their absolute values, $\mathbf{P}_2 $ is an $N_2^\partial \times N_2$ projection matrix that extracts values on boundary triangles, $\Pi_S$ is a diagonal matrix whose entries are the areas of boundary triangles, i.e., $\mathrm{diag}(S_1, S_2, \dots, S_{N_2^\partial})$, and  $\mathbb{Q}$ is an $N_2^\partial \times N_0^\partial$ that interpolates nodal values of the scalar field on the boundary to the centroids of boundary triangles, and is given by
\[
\mathbb{Q} = \frac{1}{6} \, \abs{\ded_1^\partial} \, \abs{\ded_0^\partial},
\]
where $\ded_k^\partial$ denotes the discrete exterior derivative of degree $k$ on the boundary mesh $\Gamma$, defined analogously to $\ded_k$. The $N_0^\partial\times 1$ vectors $\partial_{\nn}{\mathbf{a}}_\tau^{(\Gamma)}$ and $\partial_{\nn} {\mathbf{\Phi}}_s^{(\Gamma)}$ are the outward normal derivatives $\partial_{\nn}\tilde{A}_\tau$ and $\partial_{\nn}\Phi_s$, respectively, evaluated at the boundary nodes. Notice that the prefactor $1/6$ in Eq.~\eqref{eq:d2C} arises as the product of two terms: (1) a factor $1/2$, which ensures that the entries of $\abs{\ded_0^\mathsf{T}} \, \abs{\ded_1^\mathsf{T}}$ are either $0$ or $1$; and  
(2) a factor $1/3$, which accounts for the contribution of each primal boundary triangle to the dual area, completing the dual mesh on $\Gamma$.

Next, the column vectors $\partial_{\nn}{\mathbf{a}}_\tau^{(\Gamma)}$ and $\partial_{\nn} {\mathbf{\Phi}}_s^{(\Gamma)}$ in Eqs.~\eqref{eq:Atau} and~\eqref{eq:Phitau} must be obtained by imposing the continuity conditions of the fields at $\Gamma$, and using Eqs.~\eqref{eq:main_Ax}-\eqref{eq:main_Phi} (with $k_0^2$ scaling of Eq.~\eqref{eq:main_Phi}) to obtain 
\begin{equation}\label{eq:nA}
\partial_{\nn}{\mathbf{a}}_\tau^{(\Gamma)} = \bmsf{S}^{-1} \brc{\boldsymbol{f}_{\tilde{A}_\tau} - \bmsf{D} {\mathbf{a}}_\tau^{(\Gamma)}},\quad \forall\,\tau \in\set{x,y,z}   
\end{equation}
and
\begin{equation}\label{eq:nPhi}
\partial_{\nn}{\mathbf{\Phi}}_s^{(\Gamma)} = \bmsf{S}^{-1} \brc{  \boldsymbol{f}_{\Phi_s} - \bmsf{D} \mathbf{\Phi}_s^{(\Gamma)}}.
\end{equation}
Defining the $N_0^\Par\times N_0$ projection matrix $\mathbf{P}_0$ that maps a 0-form on $\Omega_2$ to its evaluations at the boundary nodes, we obtain
\[
 \mathbf{a}_\tau^{(\Gamma)} =  \mathbf{P}_0 \overline{\mathbf{a}}_\tau, \quad {\mathbf{\Phi}}_s^{(\Gamma)} = \mathbf{P}_0\overline{\mathbf{\Phi}}_s,
\]
where $\mathbf{P}_0$ is obtained by starting with the $N_0\times N_0$ identity matrix and deleting the rows corresponding to interior node indices. Substituting these relations into Eqs.~\eqref{eq:nA}~and~\eqref{eq:nPhi} yields
\begin{equation}\label{eq:nnA}
\partial_{\nn}{\mathbf{a}}_\tau^{(\Gamma)} = \bmsf{S}^{-1} \boldsymbol{f}_{\tilde{A}_\tau} - \bmsf{S}^{-1}\bmsf{D}\mathbf{P}_0 \overline{\mathbf{a}}_\tau
\end{equation}
and
\begin{equation}\label{eq:nnPhi}
\partial_{\nn}{\mathbf{\Phi}}_s^{(\Gamma)} =    \bmsf{S}^{-1} \boldsymbol{f}_{\Phi_s} - \bmsf{S}^{-1} \bmsf{D} \mathbf{P}_0 \overline{\mathbf{\Phi}}_s.
\end{equation}
By substituting Eqs.~\eqref{eq:nnA} and~\eqref{eq:nnPhi} into Eqs.~\eqref{eq:Atau} and~\eqref{eq:Phitau}, respectively, we obtain the system
\begin{equation}
    \label{eq:linsysA}
   \renewcommand{\arraystretch}{1.2} 
    \setlength{\arraycolsep}{6pt} 
    \begin{bmatrix}
        \fourierbb{V}                & \mathbf{0}  & \mathbf{0} & -i\hd_0(\partial_x\vare)\\
        \mathbf{0}  &  \fourierbb{V}     &  \mathbf{0} & -i\hd_0(\partial_y\vare) \\
        \mathbf{0} &  \mathbf{0}   &  \fourierbb{V}     &  -i\hd_0(\partial_z\vare)\\
        -ik_0^2\hd_0(\partial_x\vare) &  -ik_0^2\hd_0(\partial_y\vare)  & -ik_0^2\hd_0(\partial_z\vare)& \fourierbb{L}     \\
    \end{bmatrix}  \begin{bmatrix}
        \overline{\mathbf{a}}_x\\
       \overline{\mathbf{a}}_y\\
      \overline{\mathbf{a}}_z\\
    \overline{\mathbf{\Phi}}_s
    \end{bmatrix} = -\begin{bmatrix}
\ded_2^{(\mathrm{c})}\bmsf{S}^{-1}\boldsymbol{f}_{\tilde{A}_x}\\
\ded_2^{(\mathrm{c})}\bmsf{S}^{-1}\boldsymbol{f}_{\tilde{A}_y}\\
\ded_2^{(\mathrm{c})}\bmsf{S}^{-1}\boldsymbol{f}_{\tilde{A}_z}\\
\ded_2^{(\mathrm{c})}\bmsf{S}^{-1}\boldsymbol{f}_{\Phi_s}\\
    \end{bmatrix}
\end{equation}
where
\[
{\fourierbb{V}}:= -\ded_0^\mathsf{T}\hd_1\ded_0 + k_0^2\hd_0(\vare) - \ded_2^{(\mathrm{c})}\bmsf{S}^{-1}\bmsf{D}\mathbf{P}_0,
\]
and
\[
\fourierbb{L} := -\ded_0^\mathsf{T}\hd_1(\vare)\ded_0 + k_0^2\hd_0(\vare^2) - \ded_2^{(\mathrm{c})}\bmsf{S}^{-1}\bmsf{D}\mathbf{P}_0.
\]
The numerical implementation of the hybrid DEC-SIE formulation is detailed in Algorithm~\ref{alg:dec_sie_pipeline_ref}, which outlines the computational pipeline from domain discretization to the final solution.
\begin{algorithm}[H]
\caption{Hybrid DEC-SIE Algorithm}
\label{alg:dec_sie_pipeline_ref}
\begin{algorithmic}[1]
    \State \textbf{Input:} Primal mesh of $\Omega_2$, material properties $\vare(\rs)$, incident fields.
    \Comment{Step 1: Discretization of the Domain}
    \State Generate primal tetrahedral mesh of $\Omega_2$ with $N_0$ nodes.
    \State Identify boundary nodes on $\Gamma$ ($N_0^\partial$ nodes).
    \Comment{Step 2: Construction of DEC Operators}
    \State Define discrete exterior derivatives $\ded_0, \ded_1$.
    \State Define discrete Hodge stars $\hd_0, \hd_1$ based on $\vare$.
    \State Define discrete exterior derivatives on the boundary $\ded_0^\partial, \ded_1^\partial$.
    \State Construct projection matrices $\mathbf{P}_0$ and $\mathbf{P}_2$.

    \Comment{Step 3: Construction of SIE Matrices.}
    \State Discretize single-layer potential to obtain $\bmsf{S}$.
    \State Discretize double-layer potential to obtain $\bmsf{D}$.
    \State Evaluate incident field terms $\boldsymbol{f}_{\tilde{A}_\tau}$ and $\boldsymbol{f}_{\Phi_s}$.
 
    \Comment{Step 4: Linear System Formation}
    \State Compute $\fourierbb{V} = -\ded_0^\mathsf{T}\hd_1\ded_0 + k_0^2\hd_0(\vare) - \ded_2^{(\mathrm{c})}\bmsf{S}^{-1}\bmsf{D}\mathbf{P}_0$.
    \State Compute $\fourierbb{L} = -\ded_0^\mathsf{T}\hd_1(\vare)\ded_0 + k_0^2\hd_0(\vare^2) - \ded_2^{(\mathrm{c})}\bmsf{S}^{-1}\bmsf{D}\mathbf{P}_0$.
    \State Compute $\hd_0(\partial_x\vare)$, $\hd_0(\partial_y\vare)$, and $\hd_0(\partial_z\vare)$, 
    \State Assemble the global linear system~\eqref{eq:linsysA}.

    \Comment{Step 5: Solving the Linear System}
    \State Solve the linear system for $\overline{\mathbf{a}}_x, \overline{\mathbf{a}}_y, \overline{\mathbf{a}}_z, 
    \overline{\mathbf{\Phi}}_s$ using an iterative solver.
     \State \textbf{Output:} Discrete 0-forms $\overline{\mathbf{a}}_x, \overline{\mathbf{a}}_y, \overline{\mathbf{a}}_z, \overline{\mathbf{\Phi}}_s$.
\end{algorithmic}
\end{algorithm}
\section{Numerical results}\label{sec:results}
This section presents the results of a comprehensive validation study designed to assess the accuracy, stability, and robustness of the proposed hybrid DEC-SIE method for electromagnetic analysis of inhomogeneous media under diverse scenarios. The validation is performed using two distinct approaches. First, we utilize the analytical Mie series solution for a dielectric sphere embedded in free space, as well as the case of a multilayered sphere. These experiments are designed to demonstrate  
(i) convergence of the method with increasing discretization density for a fixed permittivity, (ii) robustness against high permittivity contrasts by varying the relative permittivity to large values, and (ii) broadband stability of the formulation, particularly in the low-frequency regime. For this study, we report the discrete $L^2$-norm errors over the DEC domain $\Omega_2$, defined for $u\in L^2(\Omega_2)$ as
\begin{equation}\label{eq:L2Norm}
\norm{u} := \left( \sum_{j=1}^{N_3} |u(\rs_j)|^2 \, \abs{\sigma^3_j} \right)^{1/2},
\end{equation}
where \(\abs{\sigma^3_j}\) are the volumes of 3-simplices (tetrahedra) that approximate $\Omega_2$. The relative total error compared to Mie solutions is then defined as
\[
\norm{\mathbf{E}^{\mathrm{DEC\text{-}SIE}} - \mathbf{E}^{\mathrm{Mie}} }_{rel}:= \dfrac{ \sqrt{\norm{E_x^{\mathrm{err}}}^2 + \norm{E_y^{\mathrm{err}}}^2 + \norm{E_z^{\mathrm{err}}}^2}}{ \sqrt{\norm{E_x^{\mathrm{Mie}}}^2 + \norm{E_y^{\mathrm{Mie}}}^2 + \norm{E_z^{\mathrm{Mie}}}^2} },
\]
where $\norm{E_\tau^{\mathrm{err}}} := \norm{E_\tau^{\mathrm{DEC}\text{-}\mathrm{SIE}} - E_\tau^{\mathrm{Mie}}}$ for each $\tau\in\set{x,y,z}$. In all numerical examples presented hereafter, we consider a linearly polarized incident field characterized by the polarization vector \(\boldsymbol{E}_p = \hat{x}\) and propagation direction \(\hat{\boldsymbol{u}} = -\hat{z}\). The explicit expressions for the incident fields are given by Eq.~\eqref{eq:incAPhi}. Upon assembling the linear system~\eqref{eq:linsysA}, it is solved iteratively. 
\subsection{Dielectric Sphere}
\begin{figure}[htbp!]
    \centering
    \includegraphics[width=\linewidth]{ 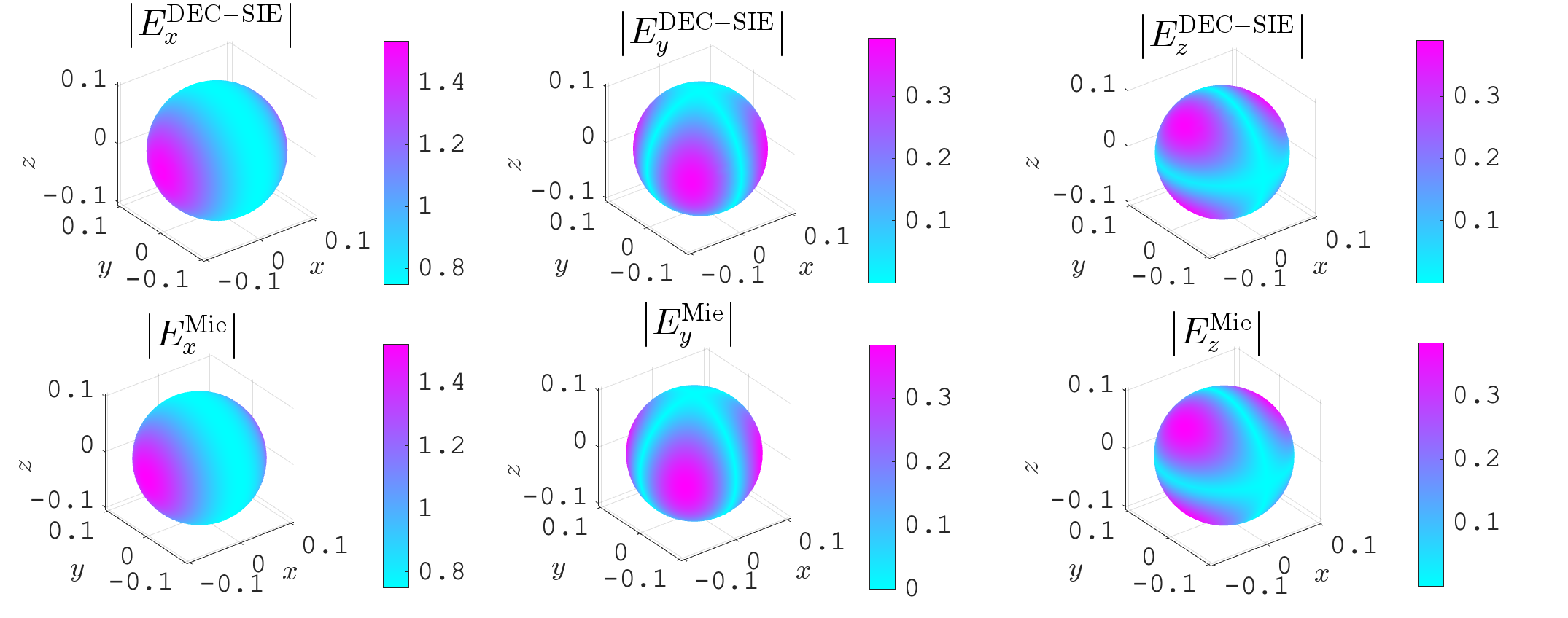}
    \caption{Comparison of the electric field component magnitudes $\abs{E_x}$, $\abs{E_y}$, and $\abs{E_z}$ on the surface $\Gamma$. The first row shows the DEC-SIE solution; the second row shows the corresponding Mie solution.}
    \label{Fig2}
\end{figure}

In this numerical experiment, we consider a dielectric sphere with relative permittivity $\vare_1 = 2.25$, and radius $a= 0.1$ m. The sphere is embedded in free space as shown in Fig.~\ref{FigT1}. The surface $\Gamma$ has to be only in free space as discussed in Section~\ref{sec2_2}. An incident plane wave with a free-space wavenumber $k_0  = 2\pi/3$ is defined by Eq.~\eqref{eq:incAPhi}. In Table~\ref{Tab1}, we report the $L^2$-norm errors of the field components $E_x$, $E_y$, and $E_z$, and relative $L^2$-norm error of the electric field $\Ev$, with mesh refinement. For this test and all subsequent tests, the condition number of the linear system matrix in Eq.~\eqref{eq:linsysA} remains consistently on the order of $10^3$, indicating good numerical conditioning of the formulation across different mesh resolutions. The field plots of the electric field component magnitudes $|E_x|$, $|E_y|$, and $|E_z|$ on the surface and the corresponding Mie solution fields are shown in Fig.~\ref{Fig2}. 
\begin{table}[htbp]
\centering
\begin{minipage}[c]{0.75\textwidth}
    \centering
    \footnotesize
    \captionof{table}{$L^2$ errors of electric field components and relative $L^2$error of $\Ev$ with mesh refinement.}
    \label{Tab1}
    \begin{tabular}{l c c c c  }
        \toprule
        $N_3$ tetrahedra & $\|E_x^{\mathrm{err}}\|$ & $\|E_y^{\mathrm{err}}\|$ & $\|E_z^{\mathrm{err}}\|$ & $\|\mathbf{E}^{\text{DEC-SIE}} - \mathbf{E}^{\text{Mie}}\|_{rel}$ \\
        \midrule
        1265   & 0.0011 & 0.0008 & 0.0011 & 0.0392 \\
        5419   & 0.0006 & 0.0005 & 0.0006 & 0.0207 \\
        13415  & 0.0004 & 0.0003 & 0.0004 & 0.0149 \\
        27924  & 0.0004 & 0.0003 & 0.0003 & 0.0117\\
        50905  & 0.0003 & 0.0002 & 0.0003 & 0.0097 \\
        81050  & 0.0003 & 0.0002 & 0.0002 & 0.0084 \\
        \bottomrule
    \end{tabular}
\end{minipage}%
\begin{minipage}[c]{0.25\textwidth}
    \centering
    \includegraphics[width=0.9\linewidth]{ 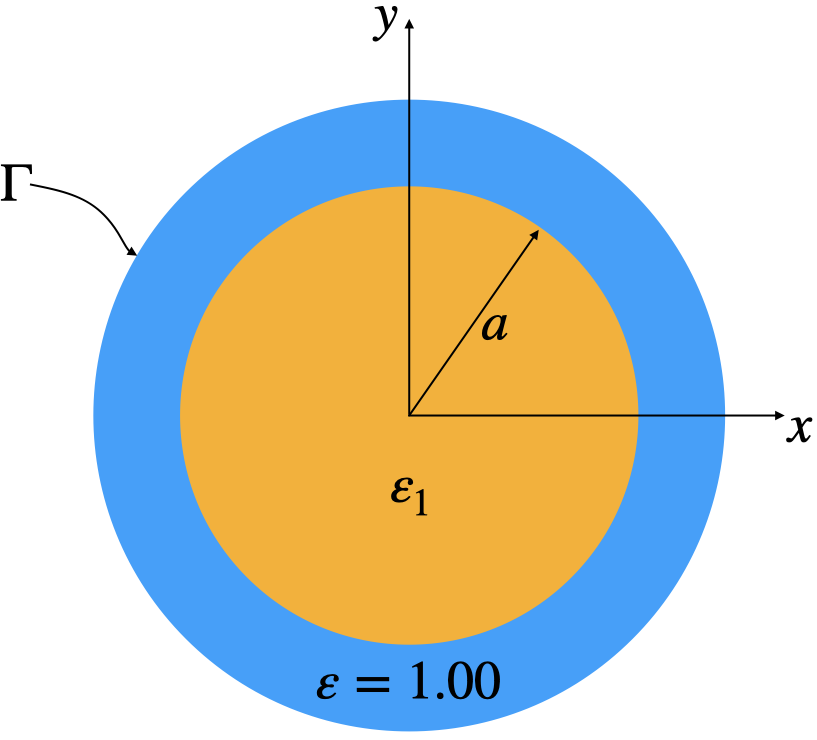}
      \footnotesize
    \captionof{figure}{Dielectric sphere.}
    \label{FigT1}
\end{minipage}
\end{table}

\begin{figure}[htbp!]
    \centering
    \includegraphics[width=\linewidth]{ 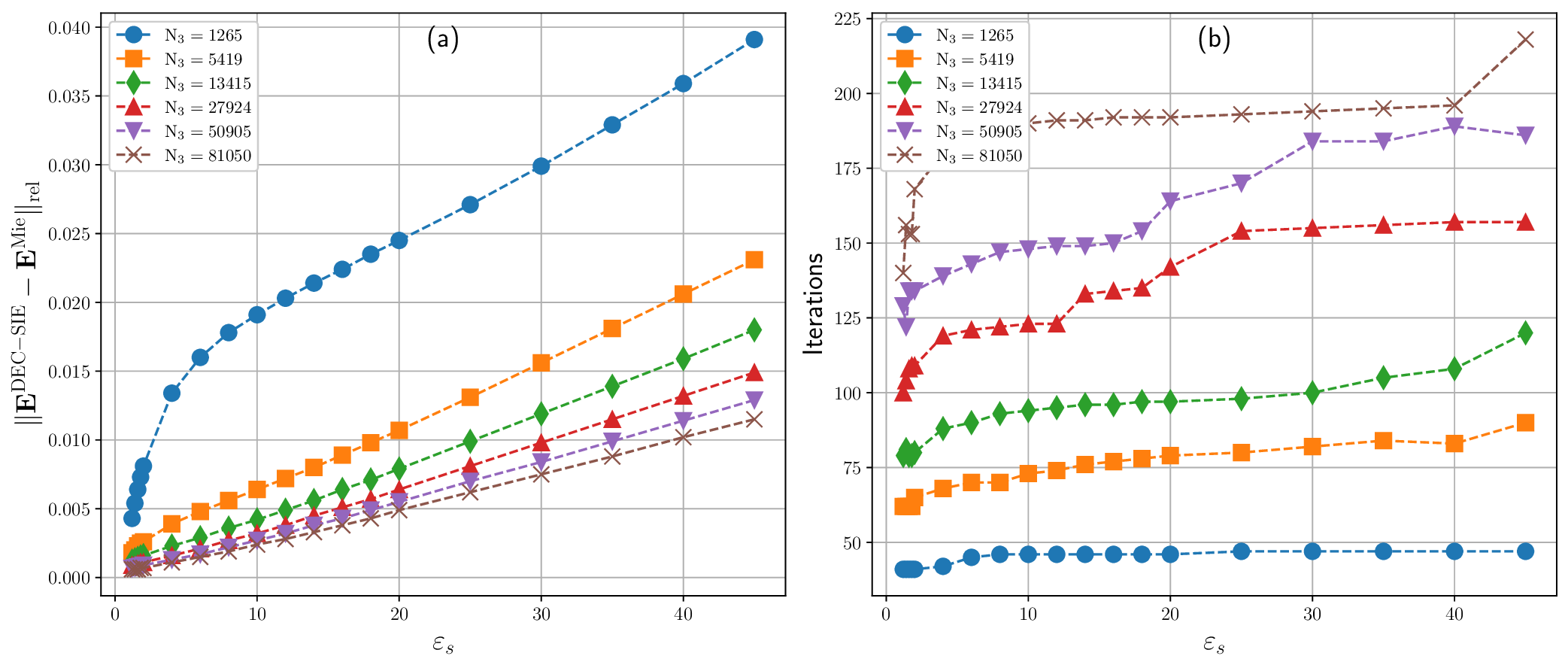}
    \caption{(a) Relative error as a function of $\vare_s \in [1, 45]$, for different mesh refinements. (b) Number of iterations required for convergence of the linear system across the same range of $\vare_s$ values.}
    \label{Fig3}
\end{figure}

For each of the meshes listed in Table~\ref{Tab1}, characterized by the number of tetrahedra $N_3$, we study the effect of high index contrast on the accuracy of the DEC-SIE method. We compute the relative error $\norm{ \mathbf{E}^{\mathrm{DEC\text{-}SIE}} - \mathbf{E}^{\mathrm{Mie}} }_{\mathrm{rel}}$ for spheres with relative permittivity $\vare_s \in [1, 45]$, under an incident field with $k_0 = \pi/25$. The resulting errors are shown as a function of $\vare_s$ in Fig.~\ref{Fig3}(a), while the corresponding number of iterations to achieve a residual relative error of approximately $10^{-13}$ of the linear system solution is plotted in Fig.~\ref{Fig3}(b). We notice the errors increase as the relative permittivity increases, while remaining within acceptable limits for the same mesh. We also notice that the errors globally decrease for all values of $\vare_s$ as the mesh is refined.

\begin{figure}[htbp!]
    \centering
    \includegraphics[width=1\linewidth]{ 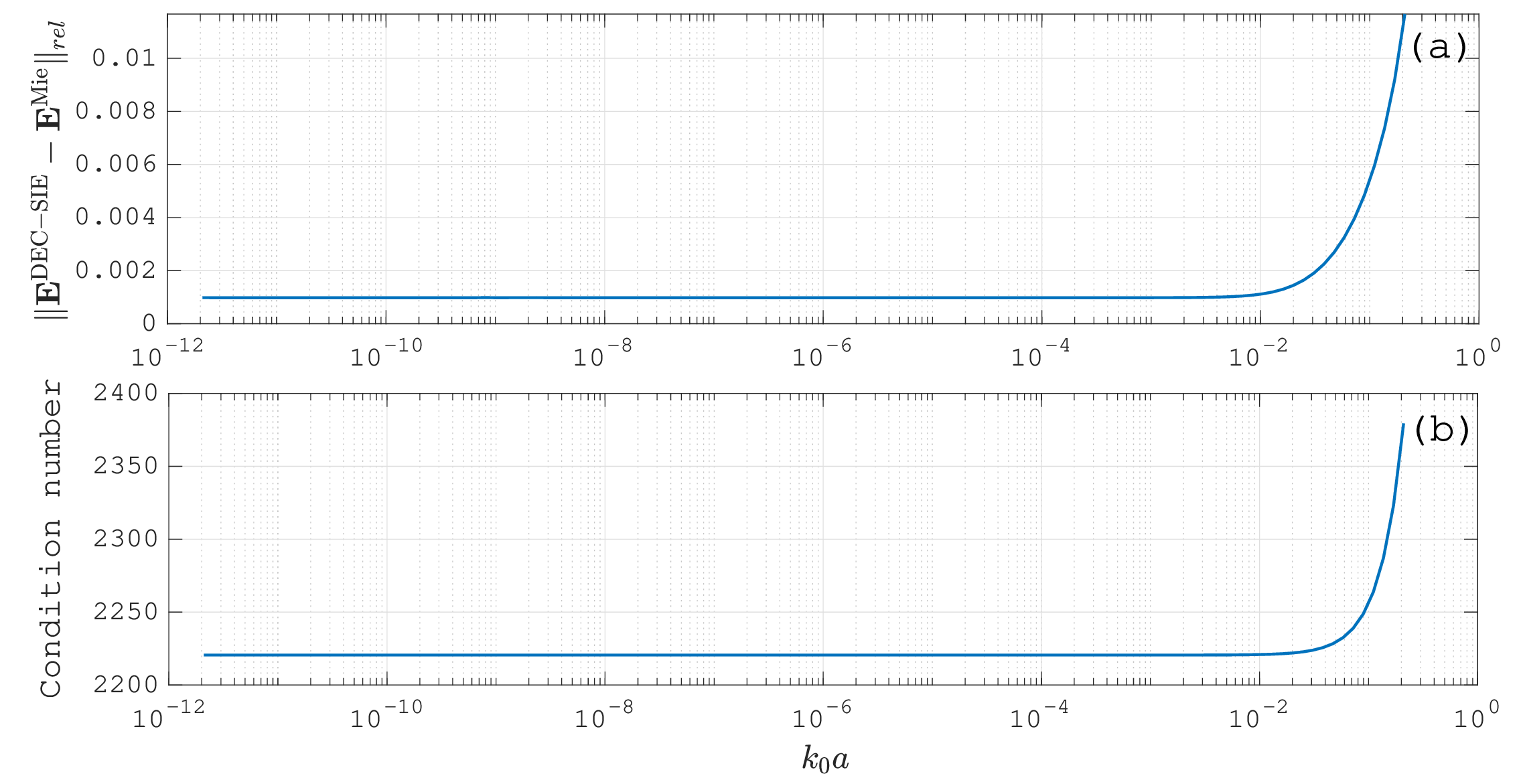}
    \caption{Broadband stability of the DEC-SIE method as $k_0 a\to 0$.  (a) Relative error as a function of the normalized frequency $k_0a$. (b) Condition number of the corresponding linear system over the same range of the $k_0a$.}
    \label{Fig4}
\end{figure}

The $\Av$-$\Phi$ formulation is immune to low-frequency breakdown. To show that the proposed potential-based DEC-SIE method is broadband stable, we consider the mesh for the dielectric sphere in Fig.~\ref{FigT1} with $N_3 = 27924$, $\vare_1 = 2.25,\, a= 0.1$ m, and compute the relative errors and the linear system condition number as $\omega \to 0$. The resulting errors are shown as a function of the normalized frequency $k_0 a$ in Fig.~\ref{Fig4}(a), while the corresponding condition number is plotted in Fig.~\ref{Fig4}(b).

\subsection{Multilayered Sphere}
In this numerical experiment, we consider a multilayered sphere with relative permittivity $\vare_1 = 2.5$, $\vare_2 = 2.25$, and radii $r_1= 0.1$ and $r_2 = 0.14$ meters, as shown in  Fig.~\ref{FigT2}. As in the previous example, the multilayered sphere is embedded in free space. An incident plane wave with a free-space wavenumber $k_0  = 2\pi/5$ is defined by Eq.~\eqref{eq:incAPhi}. In Table~\ref{Tab2}, we report the $L^2$-norm errors of the field components $E_x$, $E_y$ and $E_z$, and relative $L^2$-norm error of the electric field $\Ev$.
\begin{table}[htbp]
\centering
\begin{minipage}[c]{0.75\textwidth}
    \centering
    \footnotesize
    \captionof{table}{$L^2$ errors of electric field components and relative $L^2$error of $\Ev$ with mesh refinement.}
    \label{Tab2}
    \begin{tabular}{l c c c c}
        \toprule
        $N_3$ tetrahedra & $\|E_x^{\mathrm{err}}\|$ & $\|E_y^{\mathrm{err}}\|$ & $\|E_z^{\mathrm{err}}\|$ & $\|\mathbf{E}^{\text{DEC-SIE}} - \mathbf{E}^{\text{Mie}}\|_{rel}$ \\
        \midrule
          1549 & 0.0084 & 0.0048 & 0.0052 & 0.1378 \\
          5955 & 0.0018 & 0.0013 & 0.0014 & 0.0315 \\
        16003  & 0.0011 & 0.0008 & 0.0009 & 0.0199\\
        32929  & 0.0009 & 0.0006&  0.0007 &   0.0153\\
        58639  & 0.0007 & 0.0005 & 0.0006 & 0.0129 \\
        84917  & 0.0006 & 0.0004 & 0.0005 & 0.0115 \\
        \bottomrule
    \end{tabular}
\end{minipage}%
\begin{minipage}[c]{0.25\textwidth}
    \centering
    \includegraphics[width=0.9\linewidth]{ 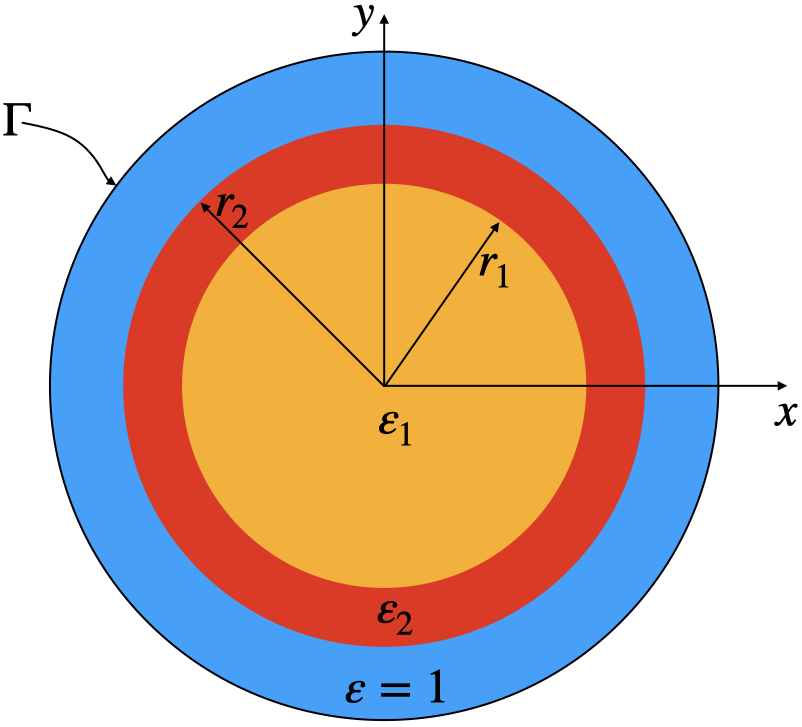}
      \footnotesize
    \captionof{figure}{Multilayered sphere.}
    \label{FigT2}
\end{minipage}
\end{table}
\subsection{A 3D scatterer}
\begin{figure}[htbp!]
    \centering
    \includegraphics[width=0.9\linewidth]{ 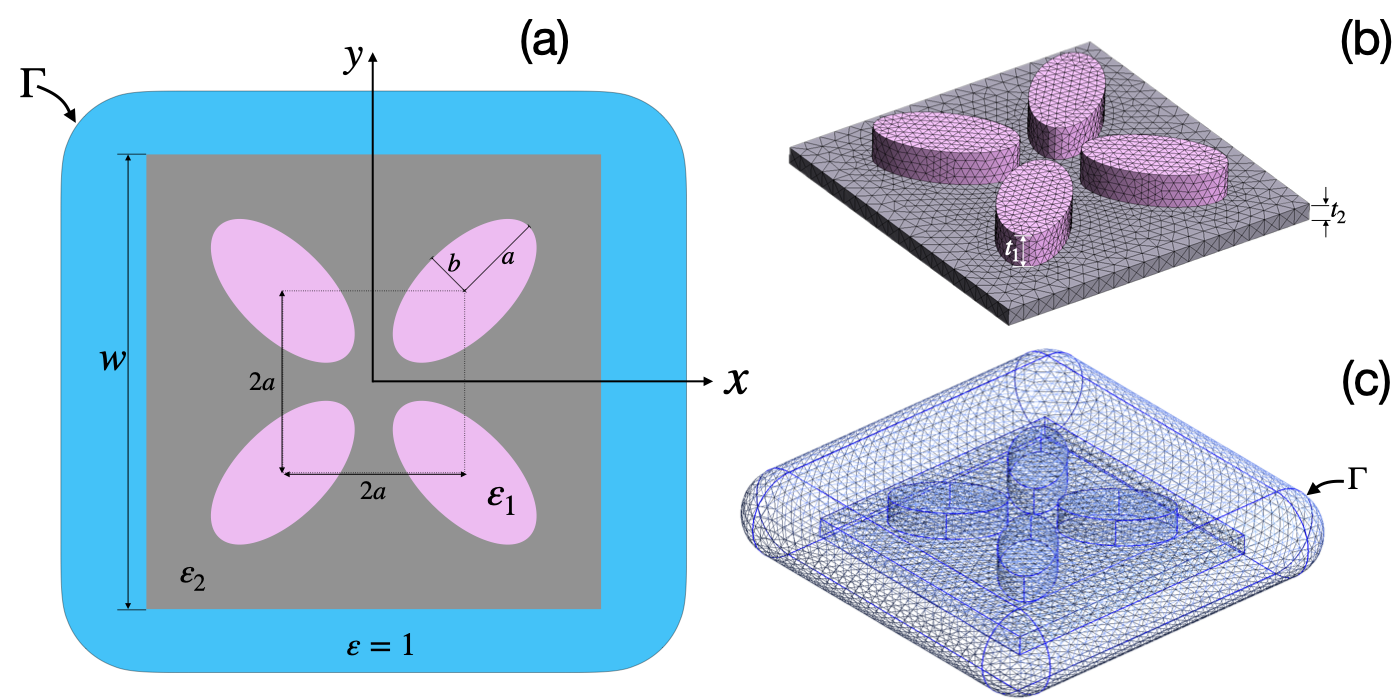}
    \caption{Geometry of the dielectric structure embedded in free space (\(\varepsilon = 1\)). \textit{(a)} Top view showing four elliptical dielectric cylinders of permittivity \(\varepsilon_1\), with semi-axes \(a\) and \(b\), arranged symmetrically with center-to-center spacing \(2a\), and placed on top of a dielectric slab of permittivity \(\varepsilon_2\). The entire structure is enclosed by the surface \(\Gamma\). \textit{(b)} Perspective view of the 3D discretized model (excluding the background mesh), illustrating the elliptical cylinders of thickness \(t_1\) situated on a dielectric slab of thickness \(t_2\), and equal width and length $w$. (c) Perspective view of the 3D discretized model showing the surface $\Gamma$ enclosing the 3D structure. }
    \label{FigT3}
\end{figure}

\begin{figure}[htbp!]
    \centering
    \includegraphics[width=\linewidth]{ 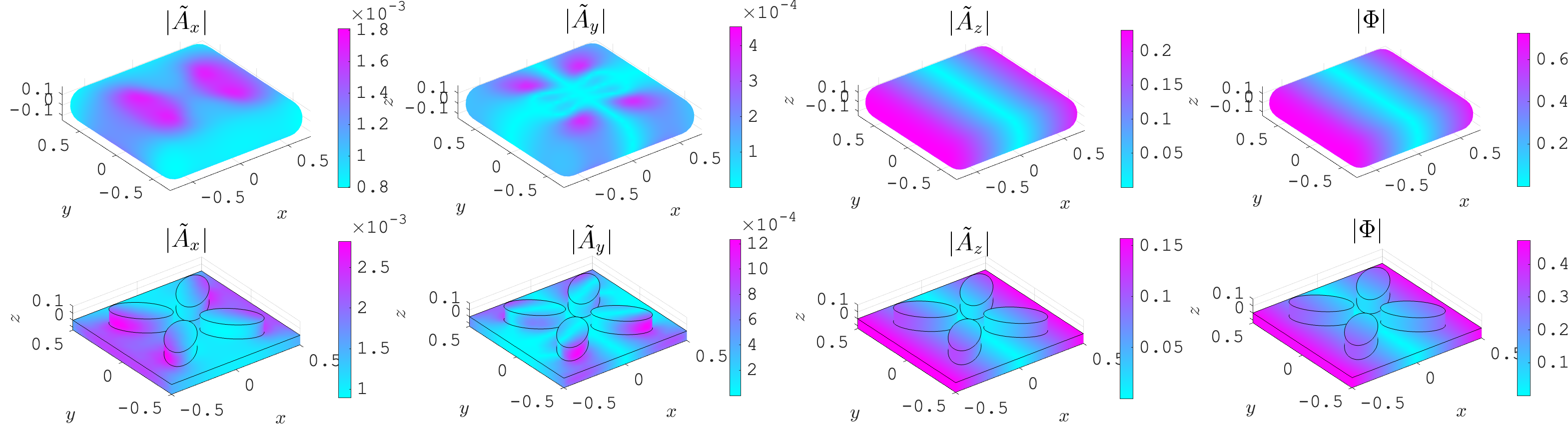}
    \caption{Surface field distributions of the computed potentials. The plots show the magnitude of each potential component, $\tilde{A}_x$, $\tilde{A}_y$, $\tilde{A}_z$, and $\Phi$. The top row corresponds to evaluations on the outer truncation boundary $\Gamma$ in free space, while the bottom row corresponds to evaluations on the surface of the interior scatterer within $\Omega_2$.}
    \label{Fig8}
\end{figure}

In this numerical experiment, we consider the 3D scatterer shown in Fig.~\ref{FigT3}. The structure consists of four elliptical dielectric cylinders of permittivity \(\varepsilon_1\), with semi-axes \(a\) and \(b\), and thickness $t_1$. The cylinders are arranged symmetrically with center-to-center spacing \(2a\), and placed on top of a dielectric slab of permittivity \(\varepsilon_2\), thickness $t_2$, and a square cross-section of side length $w$. The structure is embedded in free space, enclosed by the surface $\Gamma = \partial \Omega_2$.  
We test the hybrid DEC-SIE method with the geometry parameters (in meters) $a = 0.2$, $b = 0.1$, $t_1 = 0.1,\, t_2=0.07$, and $w = 1$. The dielectric constants are $\vare_1 =3.0$ and $\vare_2 = 2.25$. An incident plane wave with $k_0 = \pi/10$ is defined as before using Eq.~\eqref{eq:incAPhi}, and the calculated potential components $\tilde{A}_x$, $\tilde{A}_y$, $\tilde{A}_z$, and $\Phi$ are plotted on the surface $\Gamma$ as shown in the top row of  Fig.~\ref{Fig8}, while the plots for the same fields on the boundary of the internal 3D scatterer are shown in the bottom row of Fig.~\ref{Fig8}.
To assess the accuracy of the fields computed using the hybrid DEC-SIE method within the inhomogeneous domain \( \Omega_2 \), we invoke the equivalence principle and the extinction theorem~\cite{chew1995waves}. These theorems assert that the fields generated by equivalent surface potentials---obtained from surface integral equations (SIE) formulated in the exterior domain~\(\Omega_1\)---must identically cancel the incident field when evaluated inside the equivalent surface, i.e., within the interior domain~\(\Omega_2\). Any discrepancy between the directly computed interior field and the field reconstructed from boundary data, therefore, serves as a quantitative measure of numerical error.

Let \( u \in \{ \tilde{A}_x, \tilde{A}_y, \tilde{A}_z, \Phi \} \) denote a component of the computed field in \( \Omega_2 \) obtained from the DEC-SIE solution. Then, $\norm{u}$ serves as a reference magnitude of the computed field. Using the Dirichlet and Neumann traces \( u|_\Gamma \) and \( \partial_{\nn} u|_\Gamma \) on \( \Gamma \), we reconstruct the field in \( \Omega_2 \). Denote this reconstructed field by \(e_{u}\), given by
\begin{equation*}
e_{u}(\rs_j) = u^{(\mathrm{inc})}(\rs_j) - \mathcal{S}\left[ \frac{\partial u(\rs')}{\partial \hat{\boldsymbol{n}}(\rs')} \right](\rs_j) + \mathcal{D}[u(\rs')](\rs_j),
\end{equation*}
where ideally, due to the extinction theorem, we expect $e_{u}(\rs) = 0,\, \forall \, \mathbf{r} \in \Omega_2$. This field captures any residual interior response generated solely from the surface data, and thus directly reflects the consistency of the boundary data with the interior solution. The extinction relative error, $\delta_u$ is then defined by
\begin{equation*}
    \delta_u := \frac{ \norm{ e_u}   }{\norm{ u}}.
  \end{equation*}
A small value \( \delta_u \ll 1 \) indicates that the surface potentials produce negligible interior fields, validating the accuracy and consistency of the computed field \( u \) in \( \Omega_2 \). In Table 3, we report $\delta_u$ for the computed potentials as a function of $N_3$.
\begin{table}[htbp]
\centering
    \footnotesize
    \captionof{table}{Relative $L^2$ extinction errors of the potentials with mesh refinement.}
    \label{Tab3}
    \begin{tabular}{l c c c c}
        \toprule
        $N_3$  & $\delta_{\tilde{A}_x}$ & $\delta_{\tilde{A}_y}$ & $\delta_{\tilde{A}_z}$ & $\delta_{\Phi}$ \\
        \midrule
        5583 &  1.40$\times 10^{-4}$ &     2.99$\times 10^{-4}$ &     3.34$\times 10^{-5}$   &  4.31$\times 10^{-5}$ \\
        7849 &  1.36$\times 10^{-4}$ &     2.81$\times 10^{-4}$ &     3.31$\times 10^{-5}$   &  4.28$\times 10^{-5}$ \\
         9008&  1.19$\times 10^{-4}$ &     2.38$\times 10^{-4}$ &     3.22$\times 10^{-5}$   &  4.08$\times 10^{-5}$ \\
        13236 & 1.03$\times 10^{-4}$ &     2.04$\times 10^{-4}$ &     3.13$\times 10^{-5}$   &  3.91$\times 10^{-5}$ \\
        24607&  7.84$\times 10^{-5}$ &     1.54$\times 10^{-4}$ &     3.01$\times 10^{-5}$   &  3.61$\times 10^{-5}$ \\
        32132&  7.37$\times 10^{-5}$ &     1.43$\times 10^{-4}$ &     3.01$\times 10^{-5}$   &  3.56$\times 10^{-5}$ \\
        \bottomrule
    \end{tabular}
\end{table}

In the above numerical experiments, the Mie series benchmark provides an exact analytical reference for a canonical geometry with a smooth dielectric interface, enabling quantitative assessment of accuracy, convergence, and robustness to material contrast under controlled conditions. Complementarily, the extinction theorem validation rigorously verifies correct enforcement of the equivalence principle and accurate coupling across the DEC-SIE interface for general inhomogeneous structures where no analytical solution is available. The error functions directly measure the physical requirement of field cancellation inside the equivalent surface, providing a sensitive indicator of any inconsistency in the hybrid formulation. Together, these two validation frameworks comprehensively cover both analytically tractable and general geometries.
\section{Conclusion}\label{sec:conclusions}
In this work, we have presented a novel and computationally efficient hybrid DEC-SIE method for the electromagnetic analysis of heterogeneous media. The method leverages the $\Av$-$\Phi$ formulation to ensure robustness across various frequencies and material properties, offering natural compatibility with multi-physics problems and immunity to low-frequency breakdown. A key innovation of this work lies in the reformulation of the SIE, dramatically reducing the number of surface integral operators by expressing it in terms of the scalar components of the vector potential and their normal derivatives. This simplification of the SIE, combined with the adaptation of the DEC discretization to represent each Cartesian component of $\Av$ as a 0-form, enables a more efficient and consistent coupling between the DEC and SIE domains. The resulting framework provides a unified and physically consistent approach for solving complex electromagnetic problems, with significant potential for future extension to multi-physics applications, and immediate application to metasurfaces and open resonant structures.

\section*{Funding}
This work was supported in part by the National Science
Foundation under Grant No.: NSF Award CDS\&E-2202389, and in part by the Consortium on Electromagnetics Technologies (CEMT) at Purdue
University.

\appendix
\section{Exterior Calculus Representation of Maxwell's Equations}
\label{app1}

In the notation of exterior calculus, the electric field, magnetic flux density, and charge density are expressed as:
\begin{equation*}
    \begin{aligned}
        \Ef &:= E_x \ed x + E_y \ed y + E_z \ed z, \\
       \Bf &:= B_x \ed y \wedge \ed z + B_y \ed z \wedge \ed x + B_z \ed x \wedge \ed y, \\
    \boldsymbol{\rho} &:= \rho \ed x \wedge \ed y \wedge \ed z .   
    \end{aligned}
\end{equation*} 
Here $\ed$ represents the exterior derivative operator, which maps a $k$-form $\boldsymbol{\alpha}$ to a $(k+1)$-form $\ed\boldsymbol{\alpha}$. The wedge product $\wedge$ combines forms to represent oriented integrals (e.g.,  $\ed y \wedge \ed z$ defines a surface normal to the $\mathbf{x}$-axis). These forms obey antisymmetric algebraic rules $\ed x \wedge \ed x = 0$ and $\ed x \wedge \ed y = -\ed y \wedge \ed x$. A concrete example is Faraday's law, which relates a time-varying magnetic field to the electromotive force (EMF) around a closed loop.  In exterior calculus, this is written as $\ed \Ef = i\omega \Bf$. This is equivalent to writing Faraday's law $\int_{\partial \mathcal{S}} \boldsymbol{\Ef} \cdot d\mathbf{l} = i\omega\int_{\mathcal{S}} \Bf \cdot d\mathbf{s}$ , where  $\mathcal{S}$ is a surface with boundary $\partial \mathcal{S}$. Maxwell’s equations in the frequency domain become: 
\begin{equation*}
\ed \Ef = i\omega \Bf,\quad \ed \Hf = -i\omega \Df + \Jf, \quad \ed \Bf = 0, \quad \ed \Df = \tilde{\rho},
\end{equation*}
where each field is expressed as a differential form of appropriate degree. These equations hold under the condition that the domain of integration for each equation aligns with the boundary structure imposed by $\ed$. For instance, a 1-form like $\Ef$ is integrated over a curve $\partial \mathcal{S}$, which is the boundary of a surface $\mathcal{S}$ over which the 2-form $\Bf$ is defined. This structure is summarized by the generalized Stokes's theorem $ \int_{\partial \mathcal{M}} \boldsymbol{\alpha} = \int_{\mathcal{M}} \ed \boldsymbol{\alpha}$,
where $\alpha$ is a $k$-form, $\mathcal{M}$ is a $(k+1)$-dimensional oriented manifold, and $\partial \mathcal{M}$ is its boundary. Thus, the degree $k$ of a differential $k$-form determines the dimensionality of the geometric object over which it is integrated. For electromagnetic fields, the manifold $\mathcal{M} \subseteq \mathbb{R}^3$, the maximum degree of forms is limited by $k \leq 3$. 

To express electromagnetic phenomena using differential forms, we need additional operators. Here we will define the notation without specifying the precise definitions of the operators, as they ultimately are constructed within the discrete exterior calculus framework to be consistent with Maxwell's equations. The Hodge star operator $\star$, which in electromagnetism maps $k$ forms to $3-k$ forms
\[ \star: \bigwedge\nolimits^{k}(\mathcal{M}) \longrightarrow \bigwedge\nolimits^{3-k}(\mathcal{M}). \]

This complements $\ed$ and $\wedge$, enabling the definition of constitutive relations. For example, $\Ef$  and $\Hf$ are one forms and $\Df$ and $\Bf$ are two forms. The star operator in electromagnetism enables us to establish the relationship between these quantities $\Df=\star_\epsilon \Ef$ and $\Bf =\star_\mu \Hf$. Finally, the flat $\flat$  operator converts vector fields to exterior calculus 1-forms, and sharp $\sharp$ operators do the inverse.  Specifically, for vector fields $\Ev$, $\Hv$, and $\Av$, we have their corresponding 1-forms $\boldsymbol{E} = \Ev^\flat$, $\boldsymbol{H} = \Hv^\flat$, and $\boldsymbol{A} = \Av^\flat$. These additional operators provide a way to bridge exterior calculus operators with exterior derivatives. 
For a vector field $\Av$ in $\mathbb{R}^3$ and a scalar field $\Phi$, the following relations hold~\cite{mohammed2016}:
\begin{equation}\label{eq:ECI}
\begin{aligned}
(\Curl \Av)^\flat &= \star\ed\Av^\flat, \\
(\Div \Av)^\flat &= \star\ed\star\Av^\flat,\\
(\nabla \Phi)^\flat & = \ed \Phi.
\end{aligned}
\end{equation}
Each equation is a 1-form. For example, $\Av^\flat$ is a 1-form, $\ed \Av^\flat$ a 2-form, and the Hodge star on a 2-form provides a 1-form.  For readers unfamiliar with exterior calculus, references~\cite{boyuan2023,perot2014} provide a concise introduction. 

\section{The Galerkin Hodge star\label{app2}}
 
Given a tetrahedron $\mathcal{T}$ with nodes denoted locally as $[v_1,v_2,v_3,v_4]$, a Whitney 0-form $\Wf{0}_\ell$ is associated with each node $v_\ell,\, \ell = 1,\ldots, 4$, and is given by
\[
\Wf{0}_{v_\ell} = \lambda_\ell(\rs),
\]
where $\lambda(\rs)$ is the barycentric coordinates for the tetrahedron $\mathcal{T}$. Since 0-forms are scalars, there is no distinction between the Whitney 0-forms as their proxy fields. The Whitney 1-forms are associated with edges. For a tetrahedron $\mathcal{T}$, we have the locally-indexed edges $e_\ell, \,\ell = 1,\ldots,6$, where the edge $e_\ell$ is defined by the nodes forming it as $e_l = [v_{l_1},v_{l_2}]$. The Whitney 1-form associated with the edge $e_\ell$ is
\[
\Wf{1}_{e_\ell} = \lambda_{l_1}(\rs)\ed \lambda_{l_2}(\rs)-  \lambda_{l_2}(\rs)\ed \lambda_{l_1}(\rs),
\]
while its proxy field $\Wv{1}_{e_\ell}$ is
\begin{equation}\label{eq:Whitney1Proxy}
    \Wv{1}_{e_\ell} = \lambda_{l_1}(\rs)\nabla \lambda_{l_2}(\rs)-  \lambda_{l_2}(\rs)\nabla \lambda_{l_1}(\rs). 
\end{equation}
The Whitney 2-forms are associated with the triangles of $\mathcal{T}$. The triangle of $\mathcal{T}$ are $t_\ell,\, \ell = 1,\ldots,4$. As before, the triangle $t_\ell$ is defined by the nodes forming it as $t_l = [v_{l_1},v_{l_2},v_{l_3}]$ and the Whitney 2-form associated with this triangle is
\[
\Wf{2}_{t_\ell} = 2\brc{\lambda_{l_1}(\rs)\ed \lambda_{l_2}(\rs)\wedge\ed\lambda_{l_3}(\rs) +   \lambda_{l_2}(\rs)\ed \lambda_{l_3}(\rs)\wedge\ed\lambda_{l_1}(\rs) +\lambda_{l_3}(\rs)\ed \lambda_{l_1}(\rs)\wedge\ed\lambda_{l_2}(\rs)}.
\]
The proxy field $\Wv{2}_{t_\ell}$ of the 2-form $\Wf{2}_{t_\ell}$ can be written as $(\star \Wf{2})^\sharp$, and is given by~\cite{lohi2021}
\[
\Wv{2}_{t_\ell} = 2\brc{\lambda_{l_1}(\rs)\nabla \lambda_{l_2}(\rs)\times\nabla\lambda_{l_3}(\rs) +\lambda_{l_2}(\rs)\nabla \lambda_{l_3}(\rs)\times \nabla \lambda_{l_1}(\rs) +\lambda_{l_3}(\rs)\nabla  \lambda_{l_1}(\rs)\times\nabla \lambda_{l_2}(\rs)}.
\]
Given any bounded function $\xi(\rs)$ that defines material properties, e.g. $\vare(\rs)$, the Galerkin Hodge star matrices $\hd_k(\xi), \, k = 1,2$, are given by~\cite{Bossavit1999}
\begin{equation}\label{eq:ghs1}
  \brc{\hd_1(\xi)}_{i,j} = \int_{\Omega} \xi\,\Wv{1}_{e_i}\cdot\Wv{1}_{e_j} d\Omega,\quad 1 \leq i,j \leq N_1
\end{equation}
and
\begin{equation}\label{eq:ghs2}
     \brc{\hd_2(\xi)}_{i,j} = \int_{\Omega} \xi\,\Wv{2}_{t_i}\cdot\Wv{2}_{t_j} d\Omega, \quad 1 \leq i,j \leq N_2
\end{equation}
For $k = 0,3$, a diagonal representation of the operators $\hd_0$ and $\hd_3$ is used~\cite{mohamed2016}, and they are given by
\[
 \brc{\hd_0(\xi)}_{i,i} =  \abs{\ast \sigma^0_{i,\xi}},\quad 1\leq i \leq N_0
\]
\[
 \quad \brc{\hd_3(\xi)}_{i,i} = \frac{\xi_i}{\abs{\sigma^3_i}}, \quad 1\leq i \leq N_3
 \]
where $\xi_i = \xi(\rs_i)$ and $r_i$ is the centroid of tetrahedron $\mathcal{T}_i$. The diagonal entry $\abs{\ast \sigma^0_{i,\xi}}$ is the effective dual volume to the node $v_i$ accounting for the variations of the function $\xi$ over the tetrahedra, and is given by
\[
\abs{\ast \sigma^0_{i,\xi}} = \frac{1}{4}\sum_l \xi_l \abs{\sigma^3_{i,l}},
\]
where the sum runs over index $l$ of all tetrahedra sharing the node $v_i$, and the factor $1/4$ accounts for the portion of the tetrahedron volume $\abs{\sigma^3_{i,l}}$ added to the effective dual volume. 
\section{Discrete potential SIEs: method of moments}\label{app3}
Here, we briefly present the method of moments (MoM)~\cite{harrington68} used to discretize the SIEs~\eqref{eq:SIE_Ax}-\eqref{eq:SIE_Phi}. As mentioned in Section~\ref{sec:dec_sie}, the interior domain $\Omega_2$ is approximated by a simplicial complex, and the boundary of this simplicial complex---a surface triangulation---approximates $\GA$. The boundary mesh consists of $N_2^\Par$ boundary triangles, $N_1^\Par$ boundary edges, and $N_0^\Par$ boundary nodes. The position vectors of the boundary nodes are denoted by $\rs_n$, $n = 1,\ldots, N_0^\Par$. The scalar fields $\alpha$ and $\beta$ in Eq.~\eqref{eq:SIE_Ax}-\eqref{eq:SIE_Phi} are approximated by a finite expansion in nodal basis functions $h_n(\rs)$ as
\begin{equation}\label{SolApprox}
   \alpha(\rs)  \approx \sum\limits_{n = 1}^{N_0^\Par} \alpha_n h_n(\rs), \quad
   \beta(\rs)   \approx \sum\limits_{n = 1}^{N_0^\Par} \beta_n h_n(\rs),
\end{equation}
where $\alpha_n := \alpha(\rs_n)$ and $\beta_n := \beta(\rs_n)$. The basis functions $\set{h_n(\rs)}_{n = 1}^{N_0^\Par}$ are linear functions defined over each boundary triangle using barycentric coordinates, and satisfy $h_m(\rs_n) = \delta_{mn}$.

For each boundary triangle $T_m, m = 1,\ldots,N_2^\Par$, with vertex position vectors $\rs_1^{(m)}, \rs_2^{(m)}$, and $\rs_3^{(m)}$, a point $\rs\in T_m$ is associated with barycentric coordinates  $\lambda_1^{(m)}(\rs)$, $\lambda_2^{(m)}(\rs)$, and $\lambda_3^{(m)}(\rs)$. Defining
\[
\vv := \rs - \rs_3^{(m)}, \quad \vv_1 := \rs_1^{(m)}-\rs_3^{(m)}, \quad \vv_2 := \rs_2^{(m)}-\rs_3^{(m)},
\]
the barycentric functions $\lambda_1^{(m)}$, $\lambda_2^{(m)}$, and $\lambda_3^{(m)}$ are given by
\begin{equation}
    \lambda_1^{(m)}(\rs) = \frac{ \cbrc{ \vv \cdot \vv_1 } \norm{ \vv_2 }^2 - \cbrc{ \vv \cdot \vv_2 } \cbrc{ \vv_1 \cdot \vv_2 } }{ 4 S_m^2 }, \quad \rs \in T_m,
\end{equation}
\begin{equation}
    \lambda_2^{(m)}(\rs) = \frac{ \cbrc{ \vv \cdot \vv_2 } \norm{ \vv_1 }^2 - \cbrc{ \vv \cdot \vv_1 } \cbrc{ \vv_1 \cdot \vv_2 } }{ 4 S_m^2 }, \quad \rs \in T_m,
\end{equation}
and
\begin{equation}
    \lambda_3^{(m)}(\rs) = 1 - \lambda_1^{(m)}(\rs) - \lambda_2^{(m)}(\rs).
\end{equation}
Here, $S_m$ is the area of $T_m$, and we have used the identity
\begin{equation}
    \norm{ \vv_1 \times \vv_2 }^2 = \norm{ \vv_1 }^2 \norm{ \vv_2 }^2 - \cbrc{ \vv_1 \cdot \vv_2 }^2 = 4 S_m^2,
\end{equation}
which relates the triangle area to the edge vectors $\vv_1$ and $\vv_2$. 

With $h_m(\rs)$ defined, we first substitute Eqs.~\eqref{SolApprox} into Eq.~\eqref{eq:SIE_Ax}-\eqref{eq:SIE_Phi} to obtain 
\begin{equation}\label{Eq::Phi_OUT}
\frac{1}{2}\sum\limits_{n = 1}^{N_0^\Par} \alpha_n h_n(\rs) - \sum\limits_{n = 1}^{N_0^\Par} \alpha_n \brc{\int_{\GA} \frac{\partial \G (\rs,\rt)}{\partial\nn(\rt) } h_n(\rt) d\GA(\rt)} + \sum\limits_{n = 1}^{N_0^\Par} \beta_n \brc{\int_{\GA} \G (\rs,\rt) h_n(\rt) d\GA(\rt)} =  \inc{\alpha}(\rs).
\end{equation}
Second, Eq.~\eqref{Eq::Phi_OUT} is tested against $ h_m(\rs),\,m = 1,\ldots,{N_0^\Par}$, to obtain the linear system
\[
\bmsf{D}\,{\boldsymbol{\alpha}} + \bmsf{S}\,{\boldsymbol{\beta}} = \boldsymbol{f}_{\alpha}.
\]
Here, ${\boldsymbol{\alpha}} $ and ${\boldsymbol{\beta}} $ are $N_0^\Par$-dimensional column vectors, with their respective $n$th components given by $\alpha_n$ and $\beta_n$. The entries of the $N_0^\Par\times N_0^\Par$ matrices $\bmsf{D}$ and $\bmsf{S}$, denoted by $[\bmsf{D}]_{m,n}$ and $[\bmsf{S}]_{m,n}$, respectively, are given by
\begin{equation}\label{EqZ0}
    [\bmsf{D}]_{m,n} =  \frac{1}{2} \int\limits_{\GA} h_m(\rs)h_n(\rs) d\GA(\rs) -  \frac{1}{4\pi} \int\limits_{\GA} \,\int\limits_{\GA}{ h_m(\rs)h_n(\rt)   }\nn'\cdot \nabla'\frac{e^{ik_0 \abs{\rs-\rt}}}{\abs{\rs-\rt}} d\GA(\rt) d\GA(\rs),
\end{equation}
and
\begin{equation}\label{EqZ1}
   [\bmsf{S}]_{m,n} = \frac{1}{4\pi}   \int\limits_{\GA} \,\int\limits_{\GA}{ h_m(\rs)h_n(\rt)   }\frac{e^{ik_0 \abs{\rs-\rt}}}{  \abs{\rs-\rt}} d\GA(\rt) d\GA(\rs).
\end{equation}
The $m$th entry of $\boldsymbol{f}_{\alpha}$, denoted by $\brc{\boldsymbol{f}_{\alpha}}_m$ is given by
\begin{equation}\label{Eqb}
\brc{\boldsymbol{f}_{\alpha}}_m =    \int\limits_{\GA} { h_m(\rs)  \inc{\alpha}(\rs)} d\GA(\rs).
\end{equation}
For more details on the evaluation of the matrices $\bmsf{D}$ and $\bmsf{S}$, the reader may refer to~\cite{graglia1993} and the references therein. 
 



\end{document}